\documentclass{amsart}
\usepackage{amsmath,amscd,xypic,amssymb,combelow,color,enumitem,graphicx,float,comment}

\newtheorem{theorem}[subsection]{Theorem}  

\newtheorem{proposition}[subsection]{Proposition}
\newtheorem{lemma}[subsection]{Lemma}
\newtheorem{corollary}[subsection]{Corollary}

\newtheorem{definition}[subsection]{Definition}

\newtheorem{example}[subsection]{Example}
\newtheorem{remark}[subsection]{Remark}

\def\fg{{\mathfrak{g}}}

\def\fsl{{\mathfrak{sl}}}

\def\BA{{\mathbb{A}}}
\def\BC{{\mathbb{C}}}
\def\BK{{\mathbb{K}}}

\def\BN{{\mathbb{N}}}

\def\BP{{\mathbb{P}}}

\def\BQ{{\mathbb{Q}}}
\def\BZ{{\mathbb{Z}}}

\def\CL{{\mathcal{L}}}

\def\CO{{\mathcal{O}}}
\def\CP{{\mathcal{P}}}

\def\CS{{\mathcal{S}}}
\def\CT{{\mathcal{T}}}

\def\CV{{\mathcal{V}}}

\def\ph{\varphi}

\def\sym{\textrm{sym}}

\def\UU{U_q(L\fg)}

\def\UUm{U_q^-(L\fg)}

\def\bl{{\mathbf{l}}}

\def\bs{\boldsymbol{\varsigma}}

\def\nn{{\mathbb{N}^I}}
\def\zz{{\mathbb{Z}^I}}

\def\bpsi{{\boldsymbol{\psi}}}

\def\bm{{\boldsymbol{m}}}
\def\bn{{\boldsymbol{n}}}

\def\bx{{\boldsymbol{x}}}

\def\b0{{\boldsymbol{0}}}

\def\loccit{\emph{loc.~cit.~}}
\def\loccitt{\emph{loc.~cit.}}

\def\oij{\overrightarrow{ij}}

\def\bx{\boldsymbol{x}}

\def\ord{\textbf{ord }}

\def\tail{\text{tail}}
\def\head{\text{head}}
\def\nilp{\text{nilp}}

\def\etail{\emph{tail}}
\def\ehead{\emph{head}}
\def\enilp{\emph{nilp}}

\def\stab{\text{stab}}
\def\estab{\emph{stab}}

\def\gr{\text{gr}}
\def\egr{\emph{gr}}

\begin{document}

\title[Quiver moduli and quantum loop algebras]{\Large{\textbf{Quiver moduli and quantum loop algebras}}} 

\author[Andrei Negu\cb t]{Andrei Negu\cb t}

\address{École Polytechnique Fédérale de Lausanne (EPFL), Lausanne, Switzerland \newline \text{ } \ \ Simion Stoilow Institute of Mathematics (IMAR), Bucharest, Romania} 

\email{andrei.negut@gmail.com}

\maketitle
	
\begin{abstract} A classic result of Hernandez-Leclerc and Kashiwara-Kim-Oh-Park relates the $q$-characters of so-called reachable simple modules of quantum affine algebras to the Euler characteristics of certain quiver moduli spaces. We categorify and generalize this relation to all simple modules in the Hernandez-Jimbo category $\CO$ using critical $K$-theory (our results hold for quantum toroidal as well as quantum affine algebras).

\end{abstract}

\bigskip

\section{Introduction}
\label{sec:intro}

\medskip

\subsection{The setting}
\label{sub:quantum loop intro}

Let $\BK = \BQ(q)$ and fix a symmetrizable Kac-Moody Lie algebra $\fg$, with a set $I$ of simple roots. We write $(d_{ij})_{i,j \in I}$ for the symmetrized generalized Cartan matrix of $\fg$, see \eqref{eqn:cartan matrix}. The representation theory of the quantum loop algebra
$$
\UU = \BK \Big \langle e_{i,d}, f_{i,d}, \ph_{i,d'}^\pm \Big \rangle_{i \in I, d \in \BZ, d' \geq 0} \Big/ \text{relations}
$$
has long been studied, see for instance \cite{CP, HJ} for $\fg$ of finite type and \cite{N Cat} for $\fg$ of general type. We will consider the simple modules $L(\bpsi)$ of the Borel subalgebra of $\UU$ (defined in \loccitt) which are indexed by so-called highest $\ell$-weights
\begin{equation}
\label{eqn:ell weight intro}
\bpsi = \left(\psi_i(z) = \sum_{d=0}^{\infty} \frac {\psi_{i,d}}{z^d} \in \BK[[z^{-1}]]^\times \right)_{i \in I}
\end{equation}
which are rational, i.e. each $\psi_i(z)$ is the power series expansion of a rational function in $z$. Assume that the numerators and denominators of all such $\ell$-weights completely factor into $1 - \gamma z^{-1}$ for various $\gamma \in \BK^*$, which is very close to the standard treatment of these objects: typically, $\BK$ is replaced by $\BC$ and $q$ is replaced by a generic complex number. The modules $L(\bpsi)$ are completely determined by their $q$-characters (\cite{FR})
\begin{equation}
\label{eqn:q char intro}
\chi_q(V) = \sum_{\ell\text{-weights } \bpsi} \dim_{\BK} (V_{\bpsi}) [\bpsi]
\end{equation}
Above, $V_\bpsi$ denotes the generalized eigenspace of $V$ in which the commuting family of operators $\{\ph_{i,d}^+\}_{i\in I, d \geq 0}$ acts according to the eigenvalues $\{\psi_{i,d}\}_{i\in I, d\geq 0}$, and $[\bpsi]$ denotes a certain formal symbol that satisfies $[\bpsi][\bpsi'] = [\bpsi\bpsi']$ with respect to component-wise multiplication of $I$-tuples of power series. When $\fg$ is of finite type, the computation of $\chi_q(L(\bpsi))$ and its rich connections with various parts of representation theory, algebraic geometry, mathematical physics and combinatorics has been the study of much recent work (\cite{Bi, CM, FH, FM, FHOO, GHL, H, HL Cluster, KKOP, MY, Na, Q}).

\medskip

\subsection{Shuffle algebras}
\label{sub:shuffle intro}

The set $\BN$ is assumed to contain 0 throughout the present paper. Shuffle algebras have been known to provide explicit models of $\UU$ since the work of \cite{E,FO}. In more detail, there exists a subspace
$$
\CS  = \bigoplus_{\bn \in \nn} \CS_{\bn}  \subseteq \bigoplus_{\bn \in \nn} \CV_{\bn} = \CV, \quad \text{where } \CV_{\bn} = \BK [z_{i1}^{\pm 1},\dots,z_{in_i}^{\pm 1} ]^{\sym}_{i \in I}
$$
such that $\CS$ is isomorphic to the half subalgebra $\UUm \subset \UU$ generated by Drinfeld's elements $\{f_{i,d}\}_{i \in I, d \in \BZ}$. In \cite{N Cat}, we proved that the simple module $L(\bpsi)$ factors (as a graded vector space) as
\begin{equation}
\label{eqn:factor intro}
L(\bpsi) \cong L^{\text{ord }\bpsi} \otimes L^{\neq 0}(\bpsi)
\end{equation}
where the graded vector space $L^{\text{ord }\bpsi}$ only depends on the orders of the rational functions $\psi_i(z)$ at $z = 0$, and was described in \cite{N Char}. Meanwhile, the vector space $L^{\neq 0}(\bpsi)$ will be shown in \cite{HN} to match the underlying graded vector space of the same-named simple module (\cite{H Shifted}) of an appropriately shifted quantum loop algebra. We proved in \cite{N Cat} that this vector space has the following shuffle algebra description
\begin{equation}
\label{eqn:module intro}
L^{\neq 0}(\bpsi) = \CS \Big / \left( \CS \cap \Theta(\bpsi) \right)
\end{equation}
where $\Theta(\bpsi) \subseteq \CV$ is an explicit set of Laurent polynomials that we recall in Subsection \ref{sub:simple}. Moreover, the objects above admit gradings that will be denoted by
\begin{equation}
\label{eqn:grading intro}
L^{\neq 0}(\bpsi) = \mathop{\bigoplus_{\bn = (n_i \geq 0)_{i \in I}}}_{\bx = (x_{ia} \in \BK^*)^{i \in I}_{1 \leq a \leq n_i}} \CS_{\bn} \Big / \left( \CS_{\bn} \cap \Theta(\bpsi)_{\bx} \right)
\end{equation}
where we recall $\Theta(\bpsi)_{\bx} \subseteq \CV_{\bn}$ in Subsection \ref{sub:simple}. If $A_{i,x}^{-1} = \left[ \left( \frac {z-xq^{d_{ij}}}{zq^{d_{ij}}-x} \right)_{j \in I} \right]$, then
\begin{equation}
\label{eqn:second factor intro}
\chi_q(L^{\neq 0}(\bpsi)) = [\bpsi] \sum_{\bn, \bx \text{ as above}} \dim_{\BK} \left( \CS_{\bn} \Big / \left( \CS_{\bn} \cap \Theta(\bpsi)_{\bx} \right) \right)  \prod_{i \in I} \prod_{a=1}^{n_i} A_{i,x_{ia}}^{-1}
\end{equation}
see \cite[Subsection 4.8]{N Char}. Formula \eqref{eqn:factor intro} implies the equality
$$
\chi_q(L(\bpsi)) = \chi^{\text{ord }\bpsi} \cdot \chi_q(L^{\neq 0}(\bpsi))
$$
The factor $\chi^{\text{ord }\bpsi}$ was calculated in \cite{N Char} (in finite types, in accordance with a conjecture of \cite{MY, W}). Thus, the computation of $q$-characters boils down to understanding the dimensions of the vector spaces $\CS_{\bn} / \left( \CS_{\bn} \cap \Theta(\bpsi)_{\bx} \right)$. For $\fg$ of classical types, these dimensions can in principle be calculated following the residue description of \cite{N Cat}, and we expect one may recover the combinatorial formulas of \cite{J, MY Path, NN1, NN2, TDL}.

\medskip

\subsection{Euler characteristics}
\label{sub:euler intro}

The main purpose of the present paper is to produce a geometric realization of $\CS_{\bn} / \left( \CS_{\bn} \cap \Theta(\bpsi)_{\bx} \right)$, inspired by the following seminal construction of Hernandez-Leclerc. We assume that \footnote{As has long been known to experts, working with $\bpsi$ of the form \eqref{eqn:bpsi intro} does not represent any actual restriction, since the $q$-character of an arbitrary $L(\bpsi)$ naturally splits as a product of $q$-characters of $L(\bpsi)$'s of the form \eqref{eqn:bpsi intro}: in other words, poles and zeroes of $\bpsi$ that are in different geometric progressions with ratio $q$ do not interact with each other.}
\begin{equation}
\label{eqn:bpsi intro}
\bpsi = \left( \psi_i(z) = \text{constant}_i \cdot \frac {\prod_{s \in \BZ} \left (1- \frac {q^s}z\right)^{l_{i,s}}}{\prod_{s \in \BZ} \left (1- \frac {q^s}z\right)^{m_{i,s}}} \right)_{i \in I}
\end{equation}
for various $l_{i,s}, m_{i,s} \in \BN$, almost all of which are 0. We will encode these numbers in the vectors of non-negative integers $\bl = (l_{i,s})_{i \in I, s \in \BZ}$ and $\bm = (m_{i,s})_{i \in I, s \in \BZ}$.

\medskip

\begin{theorem}
\label{thm:hl}

(general conjecture in \cite{HL Cluster}, proof for Kirillov-Reshetikhin $\bpsi$ in loc. cit., general proof in \cite{KKOP}) For $\fg$ of finite type and ``reachable" highest $\ell$-weights $\bpsi$, 
\begin{equation}
\label{eqn:conj intro}
\chi_q(L(\bpsi)) = [\bpsi] \mathop{\bigoplus_{\bn = (n_i \geq 0)_{i \in I}}}_{\bx = (x_{ia} \in q^{\BZ})^{i \in I}_{1 \leq a \leq n_i}} (\text{Euler characteristic of $N^{\emph{stab}}_{\bx,\bm,\bl}$})  \prod_{i \in I} \prod_{a=1}^{n_i} A_{i,x_{ia}}^{-1}
\end{equation}
where the algebraic variety $N^{\emph{stab}}_{\bx,\bm,\bl}$ is defined as follows:

\medskip

\begin{itemize}[leftmargin=*]

\item consider the quiver $Q^{\emph{gr}}$ with vertex set $I \times \BZ$ and arrows $(i,s) \mapsto (j,s-d_{ij})$;

\medskip

\item consider the moduli space $M^{\emph{stab}}_{\bx,\bm}$ of stable $Q^{\emph{gr}}$ representations of dimension $|\{1 \leq a \leq n_i \text{ s.t. } x_{ia} = q^s\}|$ and framing $m_{i,s}$ at any vertex $(i,s)$ (see Section \ref{sec:moduli});

\medskip

\item consider the closed subvariety
\begin{equation}
\label{eqn:n in m}
N^{\emph{stab}}_{\bx,\bm,\bl} \subset M^{\emph{stab}}_{\bx,\bm}
\end{equation}
defined by intersecting the critical locus of the potential (defined below)
\begin{equation}
\label{eqn:critical intro}
\emph{Crit} \left(\emph{Tr}(W^{\emph{gr}}) \right)
\end{equation}
with the closed subvarieties defined by the vanishing condition of $l_{i,s}$ generic linear combinations of framed paths ending at any vertex $(i,s) \in I \times \BZ$  (see Section \ref{sec:moduli}).

\end{itemize}

\medskip

\noindent The potential $W^{\emph{gr}}$ in \eqref{eqn:critical intro} is the sum over all cycles of the form ($\forall i \neq j, \forall s \in \BZ$)
$$
\xymatrix{& & (j,s) \ar[rrd] & & \\
 (i,s +d_{ij}) \ar[rru] &  (i,s+d_{ij} + d_{ii}) \ar[l] & \dots \ar[l] &  (i,s-d_{ij} - d_{ii}) \ar[l] & (i,s-d_{ij}) \ar[l]}
$$
and its trace naturally induces a function on $M^{\emph{stab}}_{\bx,\bm}$. The critical locus of this function, which appears in \eqref{eqn:critical intro}, is explicitly determined by the equations
\begin{equation}
\label{eqn:w equation intro}
\left\{\frac {\partial W^{\emph{gr}}}{\partial e} = 0 \right\}_{\forall \text{ arrow }e}
\end{equation}

\end{theorem}

\medskip

\noindent Kashiwara-Kim-Oh-Park proved the result above by establishing a monoidal categorification of the cluster algebra associated to the quiver $Q^{\text{gr}}$ (for which the relation to Euler characteristics of quiver moduli had been established on quite general grounds in \cite{DWZ}), following the expectation of Hernandez-Leclerc. This explains the need for the adjective ``reachable" in the statement of Theprem \ref{thm:hl}, since these are precisely the modules whose $q$-characters correspond to cluster monomials.

\medskip

\subsection{General highest $\ell$-weights $\bpsi$}
\label{sub:toward}

We wish to generalize \eqref{eqn:conj intro} to all highest $\ell$-weights $\bpsi$, even those which produce simple modules only of the Borel subalgebra of $\UU$. Firstly, in this generality it is natural to replace $\chi_q(L(\bpsi))$ with $\chi_q(L^{\neq 0}(\bpsi))$ in the left-hand side (these two quantities are equal for highest $\ell$-weights for which the action extends to the entire $\UU$). Secondly, it is naive to hope that the right-hand side should be something as nice as the Euler characteristic. Therefore, our intended generalization will be to find a geometric formula for the coefficients of $\chi_q(L^{\neq 0}(\bpsi))$, in terms of the algebraic variety $N^{\stab}_{\bx,\bm,\bl}$. Since we already know from \eqref{eqn:second factor intro} that the aforementioned coefficients are given by the dimensions of the vector spaces $\CS_{\bn}/(\CS_{\bn} \cap \Theta(\bpsi)_{\bx})$, then the problem becomes equivalent to the following:
\begin{equation}
\label{eqn:problem 1}
\textbf{realize } \CS_{\bn} \Big/ \left(\CS_{\bn} \cap \Theta(\bpsi)_{\bx}\right) \textbf{ in terms of the geometry of } N^{\stab}_{\bx,\bm,\bl} 
\end{equation}
To set up our solution to the problem above, in what follows we fix arbitrary
$$
\bn = (n_i \geq 0)_{i \in I}  \qquad \text{and} \qquad \bx = (x_{ia} \in q^{\BZ})_{i \in I, 1 \leq a \leq n_i}
$$
Let $K(X)$ be the $0$-th (topological or algebraic) $K$-theory ring of a variety $X$. 

\medskip

\begin{theorem}
\label{thm:moduli intro}

For any symmetrizable Kac-Moody $\fg$, we have a ring isomorphism
\begin{equation}
\label{eqn:moduli intro}
\Gamma : K(M_{\bx,\bm}^{\emph{stab}}) \otimes_{\BQ} \BK \xrightarrow{\sim} \CV_{\bn} \Big/ \Theta \left( \frac 1{\bpsi^{\emph{den}}} \right)_{\bx}
\end{equation}
for any $\bn,\bx,\bm$ as above, where $\bpsi^{\emph{den}}$ denotes the denominator of \eqref{eqn:bpsi intro}. The space $M_{\bx,\bm}^{\emph{stab}}$ has no odd $K$-theory.

\end{theorem}

\medskip

\noindent Theorem \ref{thm:moduli intro} is an immediate consequence of Proposition \ref{prop:rewrite} and Theorem \ref{thm:nilp}. The latter theorem will be proved in Section \ref{sec:moduli} for moduli spaces of stable framed representations of arbitrary quivers; it complements known descriptions of the cohomology of such moduli spaces from \cite{ER} and \cite{FMo}. Let us set
\begin{equation}
\label{eqn:delta intro}
\Delta_{\bl} = \prod_{i \in I} \prod_{a=1}^{n_i} \prod_{s \in \BZ} \left(1 - \frac {q^s}{z_{ia}} \right)^{l_{i,s}} \in \CV_{\bn} 
\end{equation}
It is easy to show that multiplication by $\Delta_{\bl}$ gives an isomorphism
\begin{equation}
\label{eqn:1}
\CV_{\bn} \Big / \Theta(\bpsi)_{\bx} \xrightarrow{\sim} (\Delta_{\bl} \cdot \CV_{\bn}) \Big/ \left( (\Delta_{\bl} \cdot \CV_{\bn}) \cap  \Theta \left( \frac 1{\bpsi^{\text{den}}} \right)_{\bx} \right) \subseteq \CV_{\bn} \Big/ \Theta \left( \frac 1{\bpsi^{\text{den}}} \right)_{\bx}
\end{equation}
which clearly sends
\begin{equation}
\label{eqn:2}
 \CS_{\bn} \Big/ \left(\CS_{\bn} \cap \Theta(\bpsi)_{\bx}\right) \xrightarrow{\sim} (\Delta_{\bl} \cdot \CS_{\bn}) \Big/ \left( (\Delta_{\bl} \cdot \CS_{\bn}) \cap  \Theta \left( \frac 1{\bpsi^{\text{den}}} \right)_{\bx} \right)
\end{equation}
Therefore, the isomorphism \eqref{eqn:moduli intro} allows us to more precisely state problem \eqref{eqn:problem 1} as
\begin{multline}
\label{eqn:problem 2}
\textbf{realize } \Gamma^{-1} \left(  (\Delta_{\bl} \cdot \CS_{\bn}) \Big/ \left( (\Delta_{\bl} \cdot \CS_{\bn}) \cap  \Theta \left( \frac 1{\bpsi^{\text{den}}} \right)_{\bx} \right) \right) \\ \textbf{in terms of the geometry of } N_{\bx,\bm,\bl}^{\text{stab}} \subset M_{\bx,\bm}^{\text{stab}}
\end{multline}
For example, it is quite conceivable that for the reachable highest $\ell$-weights $\bpsi$ to which Theorem \ref{thm:hl} applies, the direct image map
\begin{equation}
\label{eqn:direct image map}
K(N_{\bx,\bm,\bl}^{\text{stab}}) \rightarrow K(M_{\bx,\bm}^{\text{stab}})
\end{equation}
is injective and its image coincides with the vector space that appears on the first line of \eqref{eqn:problem 2}. If $N_{\bx,\bm,\bl}^{\text{stab}}$ could also be shown to not have any odd $K$-theory, then the previous sentence would provide an alternate proof of Theorem \ref{thm:hl}.

\medskip 

\subsection{Critical $K$-theory}
\label{sub:k-ha intro}

Let us consider critical $K$-theory, which was developed in \cite{BFK, EP, Hi} (and more recently \cite{CTZ}) and applied to the understanding of certain simple modules $L(\bpsi)$ in \cite{VV1, VV2}; see also \cite{Pad} for the original definition of $K$-theoretic Hall algebras. We will start from a certain quiver with potential which has been studied in \cite{GLS}: let $Q$ have vertex set $I$ and a single arrow $_j\square_{i} : i \rightarrow j$ for all $i,j \in I$ (this includes a single loop at every vertex $i$). Then we consider the potential
\begin{equation}
\label{eqn:potential intro}
W = \sum_{i \neq j} {_j\square_{i}}  \underbrace{ {_i\square_{i}} \dots {_i\square_{i}}}_{-c_{ij} \text{ factors}} {_i\square_{j}}
\end{equation}
and consider the moduli stack $M_{\bn}$ of $\bn$-dimensional representations of $Q$, which is acted on by $\BC^*$ with weight $q^{d_{ij}}$ on the arrow $_j\square_i$. Since the potential $W$ is $\BC^*$-invariant, then we have a localized $\BC^*$-equivariant critical $K$-theory group
$$
K^{\BC^*}(Q,W)_{\bn} \xrightarrow{\iota_{\bn}} K^{\BC^*}(M_{\bn}) = \CV_{\bn}
$$
The map $\iota_{\bn}$ is a $\CV_{\bn}$-module homomorphism with respect to natural module structures on the domain and target. Thus, there is a well-defined notion of multiplication on $K^{\BC^*}(Q,W)_{\bn}$ by $\Delta_{\bl}$ of \eqref{eqn:delta intro}. In Section \ref{sec:k-ha}, we will encounter the diagram
\begin{equation}
\label{eqn:big diagram intro}
\xymatrix{
K^{\BC^*}(Q,W)_{\bn} \ar[d]_{\pi} \ar[r]^-{\cdot \Delta_{\bl}} & K^{\BC^*}(Q,W)_{\bn} \ar[d] \ar[r]^{\iota_{\bn}} & K^{\BC^*}(M_{\bn}) \ar[d] \\
K(Q^{\gr}, W^{\gr})^{\stab}_{\bx,\bm} \otimes_{\BQ} \BK \ar[r]^-{\cdot \Delta^{\gr}_{\bl}} & K(Q^{\gr}, W^{\gr})^{\stab}_{\bx,\bm} \otimes_{\BQ} \BK \ar[r]^-{\iota^{\gr,\stab}_{\bx,\bm}} & K(M_{\bx,\bm}^{\stab})  \otimes_{\BQ} \BK}
\end{equation}
where the spaces in the bottom row are analogues of those in the top row for the quiver of Theorem \ref{thm:hl}, intersected with the stable open locus (see Subsection \ref{sub:framed stable k-ha} for more details). In the bottom row, $\Delta_{\bl}^{\gr}$ is (a slightly renormalized version of) the product of those linear factors of \eqref{eqn:delta intro} which vanish at $z_{ia} = x_{ia}$, specifically
\begin{equation}
\label{eqn:delta gr intro}
\Delta_{\bl}^{\gr} = \prod_{i \in I} \prod_{a=1}^{n_i}  \underbrace{\left(1 - \frac 1{z_{ia}} \right)^{l_{i,s}}}_{s \text{ is defined by }x_{ia}=q^s} 
\end{equation}
The map $\pi$ in \eqref{eqn:big diagram intro} is given by restriction, while $\iota_{\bn}, \iota^{\gr,\stab}_{\bx,\bm}$ behave like push-forwards, and so diagram \eqref{eqn:big diagram intro} does not commute on the nose. However, we will show that
\begin{equation}
\label{eqn:have the same image intro}
\text{the compositions } \rightarrow \rightarrow \downarrow \text{ and } \downarrow\rightarrow\rightarrow \text{ in \eqref{eqn:big diagram intro} have the same image}
\end{equation}

\medskip

\begin{theorem}
\label{thm:main intro}

For any strongly symmetrizable \footnote{See \eqref{eqn:strongly symmetrizable} for the definition of strongly symmetrizable Kac-Moody Lie algebras. We note that all finite type and affine type Lie algebras are strongly symmetrizable, except for $\widehat{\fsl}_2$.} Kac-Moody Lie algebra $\fg$ and any highest $\ell$-weight $\bpsi$, the image of the composition $\rightarrow \rightarrow \downarrow$ in diagram \eqref{eqn:big diagram intro} is naturally isomorphic to the vector space
\begin{equation}
\label{eqn:main space}
\CS_{\bn} \Big/ \left(\CS_{\bn} \cap \Theta(\bpsi)_{\bx}\right) 
\end{equation}
and thus has dimension equal to the appropriate coefficient of $\chi_q(L^{\neq 0}(\bpsi))$. 

\end{theorem} 

\medskip 

\noindent We conjectured the following result in the first draft of the present paper, and were later informed by Cao-Okounkov-Zhou-Zhou that it follows from their general theory of stable envelopes in critical $K$-theory (\cite{COZZ}).

\medskip

\begin{theorem}
\label{conj:main}

(conjectured by the author, proved by \cite{COZZ}) The maps $\pi$ and $\iota^{\egr,\estab}_{\bx,\bm}$ in \eqref{eqn:big diagram intro} are surjective and injective, respectively.

\end{theorem}

\medskip

\noindent Combining Theorems \ref{thm:main intro} and \ref{conj:main} with statement \eqref{eqn:have the same image intro} implies the following.

\medskip

\begin{corollary}
\label{cor:main}

For strongly symmetrizable $\fg$ and any highest $\ell$-weight $\bpsi$, we have
\begin{equation}
\label{eqn:resolution}
\CS_{\bn} \Big/ \left(\CS_{\bn} \cap \Theta(\bpsi)_{\bx}\right) \cong \left[ \Delta^{\egr}_{\bl} \cdot K(Q^{\egr}, W^{\egr})^{\estab}_{\bm,\bx} \right] \otimes_{\BQ} \BK
\end{equation}

\end{corollary}

\medskip

\noindent We claim that the isomorphism \eqref{eqn:resolution} provides a solution to problems \eqref{eqn:problem 1} and \eqref{eqn:problem 2}, and thus generalizes and categorifies Theorem \ref{thm:hl} (using the fact that critical $K$-theory is categorified by derived factorization categories). To see this, note that the critical $K$-theory in the right-hand side of \eqref{eqn:resolution} is a replacement of the $K$-theory of the critical locus of $W^{\text{gr}}$, while multiplication by $\Delta^{\gr}_{\bl}$ is a derived version of imposing the vanishing of the $l_{i,s}$ linear combinations of framed paths that cut out the closed subvariety \eqref{eqn:n in m}. Thus, the right-hand side of \eqref{eqn:resolution} is a replacement of $K(N_{\bx,\bm,\bl}^{\stab}) \otimes_{\BQ} \BK$, which in turn categorifies the Euler characteristic in \eqref{eqn:conj intro}. Meanwhile, the left-hand side of \eqref{eqn:resolution} categorifies the coefficients of the $q$-character of $L(\bpsi)$, due to \eqref{eqn:grading intro}.

\medskip

\subsection{Conclusion} 

To summarize, we provide a geometric construction of all simple modules $L(\bpsi)$ in the Borel category $\CO$ of quantum loop algebras $\UU$, for any strongly symmetrizable Kac-Moody Lie algebra $\fg$ (which includes all finite and affine $\fg$ except for $\widehat{\fsl}_2$). This fits into a very active field of research, that started with the Euler characteristic construction of \cite{HL Cluster}, continued with the connections with physics explored in \cite{BZ}, and more recently involved the construction of Kirillov-Reshetikhin and prefundamental modules via critical $K$-theory in \cite{VV1,VV2}. While we consider similar objects as the latter two papers, our approach is somewhat different, and in a certain sense complementary to theirs. To explain this last comment, recall that in \loccit the authors defined a convolution product on 
\begin{equation}
\label{eqn:k-ha}
K^{\BC^*}(Q^{\gr}, W^{\gr}) = \bigoplus_{\bn \in \nn} K^{\BC^*}(Q^{\gr}, W^{\gr})_{\bn}
\end{equation}
related it to the quantum affine algebra $\UU$, and defined prefundamental modules (i.e. the case $\bl = \b0$ in \eqref{eqn:bpsi intro}) using critical $K$-theory for the same quiver with potential as us. On the other hand, our viewpoint is to start from the vector space
\begin{equation}
\label{eqn:end}
\bigoplus_{\bn,\bx}  \CS_{\bn} \Big / \left( \CS_{\bn} \cap \Theta(\bpsi)_{\bx} \right)
\end{equation}
which we already know from \cite{N Cat} yields the simple module $L(\bpsi)$, and then seek geometric realizations for the direct summands in \eqref{eqn:end}. With this in mind, our Corollary \ref{cor:main} not only provides a geometric construction of all simple modules of quantum affine algebras, but provides an explicit computation of the corresponding critical $K$-theory in terms of shuffle algebras. We note that Proposition \ref{prop:yu} also provides a shuffle algebra computation of the critical K-HA studied in \cite{VV2}.

\medskip

\subsection{Acknowledgements} I would like to thank Ben Davison, David Hernandez, Shivang Jindal, Hiraku Nakajima, Tudor Pădurariu, Markus Reineke and Yukinobu Toda for numerous important conversations on the topics of the present paper, and the authors of \cite{COZZ} for their proof of Theorem \ref{conj:main}.

\bigskip

\section{Moduli of quiver representations}
\label{sec:moduli}

\medskip

\noindent  We recall the cell decomposition of moduli spaces of stable framed quiver representations from \cite{ER}, and use it to describe the ring structure on the $K$-theory of the aforementioned moduli spaces (in a different way from the earlier work \cite{FMo}). Hereafter, we write $K(X)$ for the 0-th $K$-theory \footnote{The contents of the present Section hold equally well for topological and algebraic $K$-theory; the former is closer in spirit to \cite{HL Cluster}, while the latter is closer technically to Section \ref{sec:k-ha}.} ring of a smooth complex algebraic variety $X$, with rational coefficients. We thank Markus Reineke for sharing with us much advice pertaining to the contents of the present Section.

\medskip

\subsection{Framed quiver representations}
\label{sub:quiver}

Consider an arbitrary quiver (i.e. directed graph) $Q$ with vertex set $I$ and arrow set $E$. For any $\bn \in \nn$, the moduli stack of $\bn$-dimensional quiver representations is
\begin{equation}
\label{eqn:stack}
M_{\bn} = \left(\bigoplus_{e = \oij \in E} \text{Hom}(\BC^{n_i}, \BC^{n_j}) \right) \Big/ \prod_{i \in I} GL_{n_i}
\end{equation}
where automorphisms $(g_i)_{i \in I}$ act on  representations $(\phi_e : \BC^{n_{\tail(e)}} \rightarrow \BC^{n_{\head(e)}})_{e \in E}$ by conjugation. Given any $\bm \in \nn$, we may also consider the moduli stack of framed $\bn$-dimensional quiver representations
\begin{equation}
\label{eqn:framed stack}
M_{\bn,\bm} = \left(\bigoplus_{i \in I} \text{Hom}(\BC^{m_i}, \BC^{n_i}) \bigoplus_{e = \oij \in E} \text{Hom}(\BC^{n_i}, \BC^{n_j}) \right) \Big/ \prod_{i \in I} GL_{n_i}
\end{equation}
Thus, a framed quiver representation consists of linear maps
\begin{equation}
\label{eqn:point}
\Big(\psi_i : \BC^{m_i} \rightarrow T_i, \ \phi_e : T_{\text{tail}(e)} \rightarrow T_{\text{head}(e)} \Big)_{\forall i \in I, e \in E}
\end{equation}
where $T_i$ is a $n_i$-dimensional vector space. Such a framed quiver representation is called stable if it has no proper subrepresentation which contains the images of all the $\psi_i$'s. The open substack of stable framed quiver representations
\begin{equation}
\label{eqn:framed}
M^{\stab}_{\bn,\bm} = \left(\bigoplus_{i \in I} \text{Hom}(\BC^{m_i}, \BC^{n_i}) \bigoplus_{e = \oij \in E} \text{Hom}(\BC^{n_i}, \BC^{n_j}) \right)^{\text{stable}} \Big/ \prod_{i \in I} GL_{n_i}
\end{equation}
is actually a smooth variety. Let us consider the torus action
\begin{equation}
\label{eqn:nilpotent torus action}
\BC^* \curvearrowright M^{\stab}_{\bn,\bm}
\end{equation}
which scales all the linear maps $\phi_e$ in \eqref{eqn:point} with weigth 1. The fixed point set
\begin{equation}
\label{eqn:nilpotent}
M^{\stab, \text{nilp}}_{\bn,\bm} = \left(M^{\stab}_{\bn,\bm}\right)^{\BC^*}
\end{equation}
is well-known to be proper and isomorphic to the moduli space of stable framed nilpotent quiver representations, i.e. the set of collections \eqref{eqn:point} for which there exists $N \gg 0$ such that the composition of any $N$ of the maps $\phi_e$ is 0. 

\medskip

\subsection{Paths}
\label{sub:paths}

We will fix $\bm$ and $\bn$ throughout the present section. The moduli spaces of stable framed quiver representations were studied in \cite{ER}, where a cell decomposition that we now recall was first constructed.

\medskip

\begin{definition}
\label{def:framed path}

Whereas a path in the quiver $Q$ is a chain of arrows $e_k\dots e_1$ with $\emph{tail}(e_{s+1}) = \emph{head}(e_{s})$ for all $s \in \{1,\dots,k-1\}$, we refer to a framed path as a chain
\begin{equation}
\label{eqn:framed path}
p = e_k \dots e_1 o_{i,\ell}
\end{equation}
for some $\ell \in \{1,\dots,m_i\}$, where $i = \emph{tail}(e_1)$. Intuitively, $o_{i,\ell}$ is a symbol that corresponds to the $\ell$-th coordinate of the linear map $\psi_i$ in \eqref{eqn:point}. We will call
$$
\emph{head}(p) = \emph{head}(e_k) \in I
$$
the head of a framed path $p$ as in \eqref{eqn:framed path}.

\end{definition}

\medskip

\noindent Just like a path $e_k\dots e_1$ gives rise to a map $\phi_{e_k}\dots \phi_{e_1} : T_{\text{tail}(e_1)} \rightarrow T_{\text{head}(e_k)}$ in any quiver representation, a framed path $p$ as in \eqref{eqn:framed path} gives rise to a vector
$$
v_p = \phi_{e_k}\dots \phi_{e_1} \psi_i(\ell\text{-th standard basis vector of } \BC^{m_i}) \in  T_{\text{head}(p)}
$$
in any framed quiver representation \eqref{eqn:point}.

\medskip

\subsection{Cell decompositions}
\label{sub:cell}

Let us consider a total order on the set of framed paths $p$ as in \eqref{eqn:framed path}, which has the property that
$$
p < e p
$$
for any arrow $e$ whose tail matches the head of $p$. This order induces a lexicographic order on collections of framed paths, namely
\begin{equation}
\label{eqn:collection}
S = \Big\{p_1 < \dots < p_n \Big\}
\end{equation}
We will only consider such collections for $n = \sum_{i \in I} n_i$, given our fixed $\bn = (n_i)_{i \in I} \in \nn$. If we let $i_1,\dots,i_n \in I$ denote the heads of $p_1,\dots,p_n$ respectively, then
\begin{equation}
\label{eqn:ordering}
\bs^{i_1} + \dots + \bs^{i_n} = \bn
\end{equation}
where $\bs^i \in \nn$ denotes the $I$-tuple with entry $1$ on position $i$ and $0$ everywhere else. Whenever \eqref{eqn:ordering} holds, we call $i_1,\dots,i_n$ an ordering of $\bn$. We will only consider collections $S$ which are acceptable, in the sense that every $p_b$ is either equal to $o_{i,\ell}$ (for some $i \in I$, $\ell \in \{1,\dots,w_i\}$) or to $ep_a$ (with $a<b$ and $e \in E$ having tail equal to $\head(p_a)$). The following construction is a key result of \cite{ER} (see also \cite{FMo}).

\medskip

\begin{definition}
\label{def:cell}

For any acceptable collection $S = \{p_1<\dots<p_n\}$, consider
\begin{equation}
\label{eqn:def zs}
Z_S \subset M^{\estab}_{\bn,\bm}
\end{equation}
consisting of those stable framed quiver representations for which $v_{p_1},\dots,v_{p_n}$ form a basis of the vector spaces $\{T_i\}_{i \in I}$, and moreover there exist $\gamma_{a,b,e} \in \BC$ such that
\begin{equation}
\label{eqn:prop zs}
\phi_e(v_{p_a}) = \sum_{1\leq b \leq n | ep_a > p_b} \gamma_{a,b,e} v_{p_b}
\end{equation}
for any $a \in \{1,\dots,n\}$ and any $e \in E$ such that $\etail(e) = \ehead(p_a)$ but $ep_a \notin S$.

\end{definition}

\medskip

\noindent It was shown in \cite{ER} that the coefficients $\gamma_{a,b,e}$ in \eqref{eqn:prop zs} give rise to an isomorphism
\begin{equation}
\label{eqn:iso cell}
Z_S \cong \BA^{\# \{1 \leq a , b \leq n, e \in E | \tail(e) = \head(p_a), ep_a > p_b,  ep_a \notin S\}}
\end{equation}
and that we have a cell decomposition
\begin{equation}
\label{eqn:cell decomposition}
M^{\stab}_{\bn,\bm} = \bigsqcup_{\text{acceptable collections }S} Z_S
\end{equation}
It is well-known that the existence of the cell decomposition \eqref{eqn:cell decomposition} implies that the 0-th $K$-theory (with $\BQ$-coefficients, throughout this Section) of $M^{\stab}_{\bn,\bm}$ satisfies
\begin{equation}
\label{eqn:dim k0}
\dim_{\BQ}\left(K(M^{\stab}_{\bn,\bm})\right) = \Big| \Big\{ \text{acceptable collections} \Big\} \Big|
\end{equation}
and there is no odd $K$-theory. Recall the torus action \eqref{eqn:nilpotent torus action}, whose fixed point set is the moduli space of nilpotent stable framed quiver representations.

\medskip

\begin{proposition}
\label{prop:euler nilp}

If we let $Z_S^{\enilp} = Z_S \cap M_{\bm,\bn}^{\estab, \enilp}$, then
\begin{equation}
\label{eqn:cell decomposition nilp}
M_{\bn,\bm}^{\estab, \enilp} = \bigsqcup_{\text{acceptable collections }S} Z^{\enilp}_S
\end{equation}
is also a cell decomposition, so we conclude that
\begin{equation}
\label{eqn:dim k0 nilp}
\dim_{\BQ}\left(K(M^{\estab,\enilp}_{\bn,\bm})\right) = \Big| \Big\{ \emph{acceptable collections} \Big\} \Big|
\end{equation}
and $M^{\estab,\enilp}_{\bn,\bm}$ has no odd $K$-theory.

\end{proposition}

\medskip

\begin{proof} Recall that $Z_S^{\nilp} = Z_S^{\BC^*}$ with respect to the $\BC^*$ action that scales all the arrows with weight 1. With respect to the isomorphism \eqref{eqn:iso cell}, the induced action 
$$
\BC^* \curvearrowright \BA^{\# \{ 1 \leq a , b \leq n, e \in E | \tail(e) = \head(p_a), ep_a > p_b, ep_a \notin S\}}
$$
simply scales each coordinate $\gamma_{a,b,e}$ by a certain power of the tautological character of $\BC^*$ (the power in question is 1 plus the number of consituent arrows of the path $p_a$ minus the number of constituent arrows of the path $p_b$). Therefore, $Z_S^{\nilp}$ is isomorphic to the affine subspace determined by those coordinates for which the aforementioned power is 0. Then \eqref{eqn:dim k0 nilp} follows from \eqref{eqn:cell decomposition nilp} as \eqref{eqn:dim k0} followed from \eqref{eqn:cell decomposition}.

\end{proof}

\subsection{Interlude on quiver Grassmannians}
\label{sub:interlude 1}

The moduli space $M^{\stab}_{\bn,\bm}$ has an alternative interpretation in terms of modules over the path algebra
$$
\BC Q = \bigoplus_{\text{path }p = e_k\dots e_1} \BC p
$$
with product given by concatenation of paths (we consider the empty paths $\varnothing_i$ at each $i \in I$ to be distinct idempotents in the path algebra). There is a one-to-one correspondence between representations of the quiver $Q$ and $\BC Q$-modules. Let us consider the dual injective $\BC Q$-module with socle $i \in I$
$$
I_i^* = \bigoplus_{\text{path }p = e_k\dots e_1 \text{ starting at }i} \BC p
$$
with the action via concatenation. By a slight abuse of notation, we will write $I_i^*$ for the corresponding quiver representation, and note that it is infinite-dimensional if $Q$ has oriented cycles. In the following definition, fix a quiver representation
$$
T = \{\phi_e : T_i \rightarrow T_j\}_{e = \oij \in E}
$$

\medskip

\begin{definition}
\label{def:quiver grassmannian}

For any $\bn \in \nn$, the quiver Grassmannian is
\begin{equation}
\label{eqn:quiver grassmannian}
\emph{Gr}_{\bn}(T) = \Big\{\left(T_i \stackrel{\pi_i}\twoheadrightarrow \BC^{n_i}\right)_{i \in I} \text{ compatibly with the }\phi_e\text{'s} \Big\} \Big/ \prod_{i \in I} GL_{n_i}
\end{equation}
where the compatibility condition means that every $\phi_e$ takes $\emph{Ker } \pi_i$ to $\emph{Ker }\pi_j$ (in other words, the quotients induce a quiver representation structure on $\{\BC^{n_i}\}_{i \in I}$).

\end{definition}

\medskip

\noindent When $T$ is finite-dimensional, it is clear that $\text{Gr}_{\bn}(T)$ is a closed subset of a product of ordinary Grassmannians, and thus itself a projective variety. However, we will actually be most interested in the following setup.

\medskip

\begin{proposition}
\label{prop:quiver grassmannian is framed}

For any $\bm,\bn \in \nn$, we have an isomorphism
\begin{equation}
\label{eqn:quiver grassmannian is framed}
\emph{Gr}_{\bn} \left( \bigoplus_{k \in I} (I_k^*)^{\oplus m_k} \right) \cong M^{\estab}_{\bn,\bm}
\end{equation}

\end{proposition}

\medskip

\begin{proof} Recall that $I_k^*$ is freely generated by the action of $\BC Q$ on the empty path $\varnothing_k$ at $k$, hence $(I_k^*)_i$ has a basis indexed by paths $p : k \rightarrow i$. Thus, to give surjections
\begin{equation}
\label{eqn:surj map}
\left( \bigoplus_{k \in I} \underbrace{(I_k^*)_i \oplus \dots \oplus (I_k^*)_i}_{m_k \text{ summands}} \twoheadrightarrow \BC^{n_i} \right)_{i \in I}
\end{equation}
one needs to provide spanning sets 
\begin{equation}
\label{eqn:spanning sets}
\Big(v_p^{(k,1)},\dots,v_p^{(k,m_k)}\Big)_{\text{paths } p :k \rightarrow i} \quad \text{of} \quad \BC^{n_i}
\end{equation}
for all $i \in I$. In order for this assignment to be compatible in the sense of \eqref{eqn:quiver grassmannian}, the induced quiver maps $\phi_e : \BC^{n_{\text{head}(e)}} \rightarrow \BC^{n_{\text{tail}(e)}}$ must satisfy 
$$
\phi_e(v_p^{(k,a)}) = \begin{cases} v_{ep}^{(k,a)} &\text{if }\text{tail}(p) = \text{head}(e) \\ 0 & \text{otherwise} \end{cases}
$$
for all $k \in I$ and $a \in \{1,\dots,m_k\}$. However, the equation above means that the datum \eqref{eqn:spanning sets} is completely determined by the collection of vectors
$$
\Big(v_{\varnothing_k}^{(k,1)},\dots,v_{\varnothing_k}^{(k,m_k)}\Big) \subset \BC^{n_k}
$$
as $k \in I$ varies, which altogether generate the quiver representation $\{\BC^{n_i}\}_{i \in I}$. In turn, this datum entails the same information as a point of $M^{\stab}_{\bn,\bm}$, with $v_{\varnothing_k}^{(k,a)}$ being the image of the $a$-th standard basis vector under the framing map $\BC^{m_k} \rightarrow \BC^{n_k}$. 

\end{proof}

\medskip

\noindent If we let $\text{Gr}^{\nilp}_{\bn}(T)$ denote the closed subset of $\text{Gr}_{\bn}(T)$ such that the induced quiver representation $(\BC^{n_i})_{i \in I}$ is nilpotent, then we have
$$
\text{Gr}^{\nilp}_{\bn} \left( \bigoplus_{k \in I} (I_k^*)^{\oplus m_k} \right) \cong M_{\bn,\bm}^{\stab,\nilp}
$$
By dualizing the surjective map of quiver representations \eqref{eqn:surj map}, we see that 
\begin{equation}
	\begin{aligned}
		\text{Gr}^{\nilp}_{\bn} \left( \bigoplus_{k \in I} (I_k^*)^{\oplus m_k} \right) &= \Big\{ \bn\text{-dim nilpotent subreps of } \bigoplus_{k \in I} (I_k^{**})^{\oplus m_k} \Big\} \\ 
		&=\Big\{ \bn\text{-dim subreps of } \bigoplus_{k \in I} I_k^{\oplus m_k} \Big\}
	\end{aligned}
	\label{eqn:grassmannian nilpotent}
\end{equation}
which is the description of the quiver Grassmannian in \cite{DWZ} and \cite{HL Cluster}. 

\medskip

\subsection{Residues}
\label{sub:residues}

While acceptable collections provide a great indexing set for the cells (and thus for a basis of the $K$-theory) of $M_{\bn,\bm}^{\stab}$, we will now give a description of the latter $K$-theory's underlying ring structure. To this end, consider the ring
\begin{equation}
\label{eqn:initial v}
V_{\bn} = \BQ [z_{i1}^{\pm 1}, \dots, z_{in_i}^{\pm 1} ]^{\sym}_{i \in I}
\end{equation}
where the superscript ``sym" refers to Laurent polynomials which are symmetric in the variables $z_{i1},\dots,z_{in_i}$ for each $i$ separately (such Laurent polynomials will be called ``color-symmetric"). It is well-known that we have an isomorphism
\begin{equation}
\label{eqn:initial}
V_{\bn} \xrightarrow{\sim} K(M_{\bn}) 
\end{equation}
which sends the $k$-th elementary symmetric function in $z_{i1},\dots,z_{in_i}$ to the $k$-th exterior power of the tautological rank $n_i$ vector bundle
\begin{equation}
\label{eqn:tautological}
\xymatrix{\CT_i \ar@{.>}[d] \\
M_{\bn}}
\end{equation}
which parameterizes the vector spaces $\BC^{n_i}$ in \eqref{eqn:stack} (in what follows, we will also write $\CT_i$ for the pull-back of the above vector bundle to various variants of $M_{\bn}$, such as moduli of stable framed quiver representations). $M_{\bn}$ has no odd $K$-theory.

\medskip

\begin{definition}
\label{def:theta}

For any  $i_1,\dots, i_n \in I$, consider the set
\begin{equation}
\label{eqn:theta}
\Theta_{i_1,\dots,i_n}^{\bm} \subseteq \BQ [z_{1}^{\pm 1}, \dots, z_{n}^{\pm 1} ]
\end{equation}
consisting of Laurent polynomials $F(z_1,\dots,z_n)$ such that
\begin{equation}
\label{eqn:integral vanishing}
\underset{z_n=1}{\emph{Res}} \dots \underset{z_1=1}{\emph{Res}} \left( \frac {z_1^{d_1} \dots z_n^{d_n} F(z_1,\dots,z_n)}{\prod_{a=1}^n (z_a-1)^{m_{i_a}}} \prod_{1 \leq a < b \leq n} (z_b-z_a)^{\delta_{i_ai_b} - |\{\emph{arrows }\overrightarrow{i_ai_b}\}|} \right) = 0
\end{equation}
for all $d_1,\dots,d_n \in \BZ$. 

\end{definition}

\medskip 

\noindent It is easy to see that $\Theta_{i_1,\dots,i_n}^{\bm}$ is a finite-codimension ideal, since condition \eqref{eqn:integral vanishing} boils down to the vanishing of finitely many linear combinations of derivatives of $F$ at the point $(1,\dots,1)$. For any ordering $i_1,\dots, i_n$ of $\bn$ as in \eqref{eqn:ordering}, we have an injection
\begin{equation}
\label{eqn:above}
V_{\bn} \hookrightarrow \BQ [z_{1}^{\pm 1}, \dots, z_{n}^{\pm 1} ]
\end{equation}
given by matching each variable $\{z_a\}_{1\leq a \leq n}$ with some variable of the form $z_{i_a\bullet_a}$ in \eqref{eqn:initial v} (the choice of numbers $\bullet_a$ does not matter, as long as $\bullet_a \neq \bullet_b$ whenever $a \neq b$, $i_a = i_b$). This will allow us to intersect the ideals \eqref{eqn:theta} with the ring $V_{\bn}$, which implicitly will mean pulling back $\Theta_{i_1,\dots,i_n}^{\bm}$ under the homomorphism \eqref{eqn:above}. 

\medskip

\begin{theorem}
\label{thm:nilp}

The exterior powers of $\{\CT_i\}_{i \in I}$ induce a ring homomorphism
\begin{equation}
\label{eqn:ring hom}
\rho : V_{\bn} \rightarrow K(M_{\bn,\bm}^{\estab})
\end{equation}
The homomorphism $\rho$ is surjective and its kernel is isomorphic to
\begin{equation}
\label{eqn:kernel}
\Theta_{\bn}^{\bm} = \bigcap_{i_1,\dots,i_n \text{ ordering of }\bn} \Theta_{i_1,\dots,i_n}^{\bm} 
\end{equation}
Thus, $K \left( M_{\bn,\bm}^{\estab} \right) \cong V_{\bn}/\Theta_{\bn}^{\bm}$. Moreover, $M_{\bn,\bm}^{\estab}$ has no odd $K$-theory. 

\end{theorem}

\medskip

\noindent In particular, Theorem \ref{thm:nilp} and formula \eqref{eqn:dim k0 nilp} imply that
\begin{equation}
\label{eqn:dim quot}
\dim_{\BQ}\left(V_{\bn} / \Theta_{\bn}^{\bm} \right) = \Big| \Big\{ \text{acceptable collections} \Big\} \Big|
\end{equation}
However, logically we will prove \eqref{eqn:dim quot} and then use it to deduce Theorem \ref{thm:nilp}.

\medskip

\begin{proof} For any ordering $i_1,\dots,i_n$ of $\bn$, consider the following ``flag" version of the framed moduli space \eqref{eqn:framed}
\begin{equation}
\label{eqn:framed flag}
M^{\stab}_{i_1,\dots,i_n,\bm} = 
\end{equation}
$$
= \left(\bigoplus_{i \in I} \text{Hom}(\BC^{m_i}, \BC^{n_i}) \bigoplus_{e = \oij \in E} \text{Hom}_{i_1,\dots,i_n}(\BC^{n_i}, \BC^{n_j}) \right)^{\text{stable}} \Big/ B_{i_1,\dots,i_n}
$$
where we fix decompositions $\BC^{n_i} = \oplus^{1\leq a \leq n}_{i_a = i} \BC v_a$, and we write $\text{Hom}_{i_1,\dots,i_n}(\BC^{n_i}, \BC^{n_j})$ for the linear space of homomorphisms $\BC^{n_i} \rightarrow \BC^{n_j}$ that take each $v_a$ to a linear combination of $v_b$ with $a<b$ (similarly, $B_{i_1,\dots,i_n} \subset \prod_{i \in I} GL_{n_i}$ refers to the subgroup of automorphisms which take each $v_a$ to a linear combination of $v_b$ with $a \leq b$). In other words, \eqref{eqn:framed flag} is the moduli space of stable framed $\bn$-dimensional quiver representations which admit a full flag of subrepresentations with successive quotients being the one-dimensional representations supported at vertices $i_1,\dots,i_n \in I$ and zero maps. We have a natural map
\begin{equation}
\label{eqn:natural map}
\bigsqcup_{i_1,\dots,i_n \text{ ordering of }\bn} M^{\stab}_{i_1,\dots,i_n,\bm} \xrightarrow{\upsilon} M^{\stab}_{\bn,\bm}
\end{equation}
whose image is $M^{\stab,\nilp}_{\bn,\bm}$ (this is obvious from the quiver moduli description: a nilpotent quiver representation is precisely one which can be filtered such that the successive quotients are one-dimensional quiver representations with zero maps). 

\medskip

\begin{lemma} 
\label{lem:tower} 

We have an isomorphism
\begin{equation}
\label{eqn:tower}
K (M^{\emph{stab}}_{i_1,\dots,i_n,\bm}) \cong \BQ [ z_1^{\pm 1},\dots,z_n^{\pm 1} ] \Big/ \Theta^{\bm}_{i_1,\dots,i_n}
\end{equation}
and $M^{\emph{stab}}_{i_1,\dots,i_n,\bm}$ has no odd $K$-theory. 

\end{lemma} 

\medskip

\noindent We will prove Lemma \ref{lem:tower} as soon as we finish the proof of Theorem \ref{thm:nilp}. To the latter end, consider the maps
$$
V_{\bn} \xrightarrow{\rho} K(M_{\bn,\bm}^{\stab}) \xrightarrow {\upsilon^*} \mathop{\bigoplus_{i_1,\dots,i_n}}_{\text{ordering of }\bn} K(M^{\stab}_{i_1,\dots,i_n,\bm}) \cong \mathop{\bigoplus_{i_1,\dots,i_n}}_{\text{ordering of }\bn} \BQ [ z_1^{\pm 1},\dots,z_n^{\pm 1} ] \Big/ \Theta^{\bm}_{i_1,\dots,i_n}
$$
whose composition corresponds to sending a symmetric polynomial in $\{z_{i1},\dots,z_{in_i}\}_{i \in I}$ to the same symmetric polynomial in $\{z_1,\dots,z_n\}$ (with each $z_a$ identified with $z_{i_a\bullet_a}$ for some numbers $\bullet_a$, whose particular value does not matter as long as $\bullet_a \neq \bullet_b$ whenever $a \neq b, i_a = i_b$). Therefore, we conclude that 
$$
\text{Ker }\rho \subseteq \bigcap_{i_1,\dots,i_n \text{ ordering of }\bn}  \Theta^{\bm}_{i_1,\dots,i_n} = \Theta_{\bn}^{\bm} 
$$
and so 
$$
\dim_{\BQ} \left( K(M^{\stab}_{\bn,\bm}) \right) \geq \dim_{\BQ}  \left(V_{\bn} / \text{Ker }\rho \right) \geq \dim_{\BQ} \left(V_{\bn} / \Theta_{\bn}^{\bm} \right)
$$
Once we prove \eqref{eqn:dim quot}, then formula \eqref{eqn:dim k0 nilp} will imply that all of the inequalities above are actually equalities, and the proof of Theorem \ref{thm:nilp} will be complete. 

\medskip

\noindent In order to give a direct proof of \eqref{eqn:dim quot}, let us consider the following linear map for any acceptable collection $S = \{p_1 < \dots < p_n\}$ with $i_a = \head(p_a)$:
$$
\mu_S : V_{\bn} \rightarrow \BQ, \qquad \mu_S(F) = 
$$
\begin{equation}
\label{eqn:mu}
\underset{z_n=1}{\text{Res}} \dots \underset{z_1=1}{\text{Res}} \left( F(z_1,\dots,z_n) \frac { \prod_{1 \leq a < b \leq n} (z_b-z_a)^{\delta_{i_ai_b} - | \{\text{arrows } e = \overrightarrow{i_ai_b} | ep_a \leq p_b \}|}}{\prod_{a=1}^n (z_a-1)^{| \{1 \leq \ell \leq m_{i_a} | o_{i_a,\ell} \leq p_a\}|}} \right) 
\end{equation}
Recall the lexicographic total order on acceptable collections, and let
$$
V_{\bn, \leq S} = \bigcap_{S' \leq S} \text{Ker }\mu_{S'} \qquad V_{\bn, < S} = \bigcap_{S' < S} \text{Ker }\mu_{S'}
$$
Formula \eqref{eqn:dim quot} is then an immediate consequence of the following

\medskip

\noindent \textbf{Claim 1:} $\Theta_{\bn}^{\bm} = \bigcap_{S \text{ acceptable collection}} \text{Ker }\mu_S$ 

\medskip

\noindent \textbf{Claim 2:} $\dim_{\BQ} \left(V_{\bn, < S} \Big/ V_{\bn,\leq S} \right) = 1$ for any acceptable collection $S$.

\medskip

\noindent Let us prove Claim 1: the inclusion $\subseteq$ is trivial, since the poles that arise in \eqref{eqn:integral vanishing} are at least as bad as the poles that arise in \eqref{eqn:mu}. However, when we take in \eqref{eqn:mu} the maximal acceptable collection $S = \{p_1<\dots<p_n\}$ with $\head(p_a) = i_a$, the poles in question are exactly the same as the poles that arise in \eqref{eqn:integral vanishing}. Therefore, the inclusion $\supseteq$ is an immediate consequence of the fact that 
\begin{equation}
\label{eqn:implication}
\mu_{S'}(F) = 0 \quad \forall S' \leq S \quad \Rightarrow \quad \mu_{S'}(FG) = 0 \quad \forall S' \leq S
\end{equation}
for all Laurent polynomials $F,G$ (in other words, $V_{\bn,\leq S}$ is an ideal). By induction on $S'$, it suffices to show that $\mu_S(FG) = 0$ for all Laurent polynomials $G$. By induction on $b \in \{1,\dots,n\}$, it suffices to consider (by Taylor expansion at $(1,\dots,1)$) either

\medskip

\begin{itemize}

\item  $G = z_b - 1$ if $p_b = o_{i_b,\ell}$ for some $\ell$, or

\medskip

\item $G = z_b-z_a$ if $p_b = ep_a$ for some $a<b$ and arrow $e = \overrightarrow{i_ai_b}$. 

\medskip

\end{itemize}

\noindent In either of these cases, we have
$$
\mu_{S}(FG) = \mu_{S'}(FH)
$$
for some Laurent polynomial $H$, where $S'$ is defined as the maximal acceptable collection with the same first $b-1$ paths as $S$, but whose $b$-th path is the maximal predecessor of $p_b$ in the set $(\{o_{i_b,\ell}\}_{1 \leq \ell \leq m_{i_b}} \cup \{ep_a\}_{a < b, e: i_a \rightarrow i_b} ) \backslash S$. Since $S' < S$, this completes the proof of \eqref{eqn:implication}.

\medskip

\noindent In order to prove Claim 2, it is clear that the dimension in question is $\leq 1$, because any element of $V_{\bn, < S} / V_{\bn,\leq S}$ is uniquely determined by its image under $\mu_S$. To prove equality, it suffices to produce $F \in V_{\bn}$ which is annihilated by all $\{\mu_{S'}\}_{S' < S}$ but is not annihilated by $\mu_S$. To this end, consider
$$
F_S = \text{Sym} \left[ \frac {\prod_{a=1}^n (z_a-1)^{| \{1 \leq \ell \leq m_{i_a} | o_{i_a,\ell} < p_a\}|} \prod_{1 \leq a < b \leq n} (z_b-z_a)^{|\{e = \overrightarrow{i_ai_b} | ep_a < p_b \}|}}{\prod_{1\leq a < b \leq n, \ i_a = i_b} (z_b-z_a)} \right]
$$
where ``Sym" refers to averaging over all ways to permute pairs of variables $z_a$ and $z_b$ for which $i_a = i_b$. When evaluating $\mu_{S'}(F_S)$ for any $S' = \{p_1' < \dots < p_n' \}$ using formula \eqref{eqn:mu}, the answer is equal to the sum of
$$
\underset{z_n=1}{\text{Res}} \dots \underset{z_1=1}{\text{Res}} \frac { \prod_{1 \leq a < b \leq n} (z_b-z_a)^{|\{\text{arrows } e | e\min(p_{\sigma(a)}, p_{\sigma(b)}) < \max(p_{\sigma(a)}, p_{\sigma(b)}) \}| - |\{ e = \overrightarrow{i_ai_b} | ep_a' \leq p_b' \}|}}{\prod_{a=1}^n (z_a-1)^{| \{ 1 \leq \ell \leq m_{i'_a} | o_{i_a',\ell} \leq p'_a \} | - | \{ \ell \in \{1,\dots,m_{i_a}\} | o_{i_a',\ell} < p_{\sigma(a)} \} |}} 
$$
as $\sigma \in S_n$ goes over all permutations such that $i_{\sigma(a)} = i'_a$ for all $a \in \{1,\dots,n\}$ (we write $i_a' = \head(p_a')$). It is clear that the residue above can be non-zero only if

\medskip

\begin{itemize}
	
\item $p_1' > p_{\sigma(1)}$, or

\medskip

\item $p_1' = p_{\sigma(1)}$ and $p_2' > p_{\sigma(2)}$, or

\medskip

\dots

\medskip

\item $p_1' = p_{\sigma(1)}$, $p_2' = p_{\sigma(2)}$, \dots, $p_{n-1}' > p_{\sigma(n-1)}$, or

\medskip

\item $p_1' = p_{\sigma(1)}$, $p_2' = p_{\sigma(2)}$, \dots, $p_{n-1}' = p_{\sigma(n-1)}$, $p_n' \geq  p_{\sigma(n)}$.

\end{itemize} 

\medskip 

\noindent If $S' < S$, then the inequalities above are impossible, so $\mu_{S'}(F_S) = 0$. Meanwhile if $S' = S$, the inequalities above show that only the identity permutation $\sigma \in S_n$ produces a non-zero residue (whose value is $\pm 1$), and therefore $\mu_S(F_S) \neq 0$. \end{proof}

\medskip

\begin{proof} \emph{of Lemma \ref{lem:tower}:} If $\CT$ is a rank $r$ vector bundle over a space $X$, we may define the projective space bundle $\BP_X\CT = \text{Proj}_X(\text{Sym}^\bullet(\CT)) \rightarrow X$ which parameterizes one-dimensional quotients of $\CT$. It is well-known that
$$
K(\BP_X\CT) = K(X)[z^{\pm 1}] \Big/ \wedge^\bullet\left(\frac {\CT}z\right)
$$
where the variable $z$ has $K$-theoretic degree 0, and we write
\begin{equation}
\label{eqn:wedge powers}
\wedge^\bullet\left(\frac {\CT}z\right) = \sum_{i=0}^r \frac {[\wedge^i\CT]}{(-z)^i} \qquad \wedge^\bullet\left(\frac z{\CT}\right) = \sum_{i=0}^r (-z)^i [\wedge^i\CT^\vee]
\end{equation}
In other words, if we formally write $\CT = \alpha_1 + \dots + \alpha_r$, then we have 
$$
K(\BP_X\CT) = K(X)[z^{\pm 1}] \Big / \left \{f(z) \Big| \sum_{k=1}^r \underset{z = \alpha_k}{\text{Res}} \left( \frac {z^df(z)}{\wedge^\bullet(\CT/z)} \right) = 0, \forall d \in \BZ \right\}
$$
With this in mind, formula \eqref{eqn:tower} reduces (via induction on $n$) to the claim that the forgetful map
$$
M^{\stab}_{i_1,\dots,i_n,\bm} \rightarrow M^{\stab}_{i_1,\dots,i_{n-1},\bm}
$$
is the projectivization of a vector bundle with $K$-theory class equal to
$$
[\CO^{m_{i_a}}] + \sum_{a<n} [\CL_a]^{|\text{arrows  }\overrightarrow{i_ai_n}|} - \sum_{a<n} [\CL_a]^{\delta_{i_ai_n}}
$$
where $\CL_1,\dots,\CL_{n-1}$ are the line bundles on $M^{\stab}_{i_1,\dots,i_{n-1},\bm}$ that parameterize the lines spanned by the vectors $v_1,\dots,v_{n-1}$, respectively. The aforementioned claim is straightforward, and left to the reader. \end{proof} 

\medskip

\section{$q$-characters via shuffle algebras}
\label{sec:shuffle}

\medskip

\noindent We recall the shuffle algebra incarnation of quantum loop algebras $\UU$ and of the simple modules $L(\bpsi)$ in the Borel category $\CO$. The upshot will be a formula for the $q$-character $\chi_q(L^{\neq 0}(\bpsi))$ as the dimension of an explicit vector space of Laurent polynomials. In order to better compare the contents of the present Section with the more geometric Section \ref{sec:k-ha}, we work over the ground field $\BK = \BQ(q)$, and consider only simple modules $\UU \curvearrowright L(\bpsi)$ with highest $\ell$-weight $\bpsi$ of the form \eqref{eqn:bpsi}, i.e. the numerator and denominator of $\bpsi$ completely factor into terms of the form $\{1-\gamma z^{-1}\}_{\gamma \in \BK^*}$. This is equivalent to the standard treatment of the Borel category $\CO$, in which $\BK$ is replaced by $\BC$ and $q$ is replaced by a generic complex number.

\subsection{Basic definitions}
\label{sub:basic}

Fix a finite set $I$ and a generalized Cartan matrix
\begin{equation}
	\label{eqn:cartan matrix}
C=	\left(c_{ij} = \frac {2d_{ij}}{d_{ii}} \in \BZ\right)_{i,j \in I}
\end{equation}
where $d_{ii} = 2d_i$ are positive integers with greatest common divisor 2, and $d_{ij} = d_{ji}$ are non-positive integers $\forall i \neq j$. This data corresponds to a symmetrizable Kac-Moody Lie algebra $\fg$, which we will refer to as the ``type" of all algebras considered hereafter. We interpret $I$ as a set of simple roots of $\fg$, hence we interpret $\zz$ as the root lattice. We call $\fg$ finite type if the symmetric matrix $\left(d_{ij} \right)_{i, j \in I}$ is positive definite, since in this case $\fg$ is a finite-dimensional semisimple Lie algebra. Similarly, $\fg$ is called simply-laced if the generalized Cartan matrix \eqref{eqn:cartan matrix} satisfies
\begin{equation}
\label{eqn:symmetric matrix}
d_{i} = 1, \quad \ \forall i \in I
\end{equation}
Following the suggestion of David Hernandez, we call $\fg$ strongly symmetrizable if 
\begin{equation}
\label{eqn:strongly symmetrizable}
d_{ij} \in \{ 0, - \max(d_i,d_j) \}
\end{equation}
for all $i\neq j$. All finite and affine type $\fg$ are strongly symmetrizable, except for $\widehat{\fsl}_2$.

\subsection{The big shuffle algebra}
\label{sub:big shuffle}

We now review the trigonometric version (\cite{E}) of the Feigin-Odesskii shuffle algebra (\cite{FO}) associated to the Kac-Moody Lie algebra $\fg$. Consider the vector space of Laurent polynomials in arbitrarily many variables
\begin{equation}
\label{eqn:big shuffle}
\CV = \bigoplus_{\bn \in \nn} \CV_{\bn}, \quad \text{where} \quad \CV_{(n_i \geq 0)_{i \in I}} = \BK[z_{i1}^{\pm 1},\dots,z_{in_i}^{\pm 1}]^{\text{sym}}_{i \in I} 
\end{equation}
Above, ``sym" refers to color-symmetric Laurent polynomials, meaning that they are symmetric in the variables $z_{i1},\dots,z_{in_i}$ for each $i \in I$ separately (the terminology is inspired by the fact that $i \in I$ is called the color of the variable $z_{ia}$). We make the vector space $\CV$ into a $\BK$-algebra via the following shuffle product:
\begin{equation}
	\label{eqn:mult}
	F( z_{i1}, \dots, z_{i n_i})_{i \in I} * F'(z_{i1}, \dots,z_{i n'_i})_{i \in I} = \frac 1{\bn!  \bn'!}\,\cdot
\end{equation}
$$
\textrm{Sym} \left[ F(z_{i1}, \dots, z_{in_i}) F'(z_{i,n_i+1}, \dots, z_{i,n_i+n'_i})
\prod_{i,j \in I} \mathop{\prod_{1 \leq a \leq n_i}}_{n_j < b \leq n_j+n_j'} \zeta_{ji} \left( \frac {z_{jb}}{z_{ia}} \right) \right]
$$
where we fix a total order $<$ on $I$ and define \footnote{Note that the function $\zeta_{ij}$ below is equal to $1-x^{\delta_{i < j} - \delta_{i > j}}$ times the same-named function studied in \cite{N Cat}. This means that the algebra $\CV$ of \eqref{eqn:big shuffle} is isomorphic to the same-named algebra in \loccitt, with the isomorphism given by multiplication with $\prod_{i < j} \prod_{a = 1}^{n_i} \prod_{b=1}^{n_j} \left(1 - \frac {z_{ia}}{z_{jb}} \right)$.}

\begin{equation}
\label{eqn:def zeta}
\zeta_{ij}(x) = \begin{cases} q^{-d_{ij}} - x &\text{if } i< j \\ \frac {q^{-d_{ij}}-x}{1-x} &\text{if }i = j \\ 1 - q^{-d_{ij}} x^{-1} &\text{if } i > j \end{cases}
\end{equation}
The word ``Sym" in \eqref{eqn:mult} denotes symmetrization with respect to the
\begin{equation*}
	(\bn+\bn')! := \prod_{i\in I} (n_i+n'_i)!
\end{equation*}
permutations of the variables $\{z_{i1}, \dots, z_{i,n_i+n'_i}\}$ for each $i$ independently. The reason for formula \eqref{eqn:mult} is to ensure that there exists an algebra homomorphism
\begin{equation}
\label{eqn:upsilon}
\Upsilon : \UUm \rightarrow \CV
\end{equation}
whose domain is the negative half (in Drinfeld's new realization, \cite{Dr}) of the quantum loop algebra $\UU$ of $\fg$. We refer to \cite{N Cat} for an overview of the generators-and-relations presentation of $\UU$, and note that (by definition) it involves sufficiently many relations so that $\Upsilon$ is injective. The corresponding set of relations is so far only known explicitly for finite type (\cite{B, Dr}), more generally for strongly symmetrizable (\cite{N new new}), and for simply-laced (\cite{N Loop}) symmetrizable Kac-Moody Lie algebras $\fg$. See \cite{N Arbitrary} for a discussion.

\medskip

\subsection{The (small) shuffle algebra}
\label{sub:small shuffle}

We will refer to
\begin{equation}
\label{eqn:spherical}
\CS = \text{Im }\Upsilon
\end{equation}
as simply the shuffle algebra, i.e. the subalgebra of $\CV$ generated by $\{z_{i1}^d \in \CV_{\bs^i}\}_{i \in I, d\in \BZ}$. Note that $\CS$ was denoted by $\CS^-$ in \cite{N Cat, N Char}, in order to differentiate it from its opposite algebra $\CS^{\text{op}} = \CS^+$; this distinction will not play any role in the present paper. Explicitly, $\CS$ is the $\BK$-linear span of Laurent polynomials of the form
\begin{equation}
\label{eqn:generators}
\text{Sym} \left[ z_1^{d_1} \dots z_n^{d_n} \prod_{1 \leq a < b \leq n}\zeta_{i_bi_a} \left( \frac {z_b}{z_a} \right) \right]
\end{equation}
as $i_1,\dots,i_n$ run over orderings of $\bn$, and $d_1,\dots,d_n$ run over $\BZ$ (in the above expression, we identify each variable $z_a$ with some variable of the form $z_{i_a\bullet_a}$ in the notation of \eqref{eqn:big shuffle}; the values of the numbers $\bullet_a$ do not matter due to the presence of the symmetrization, as long as we have $\bullet_a \neq \bullet_b$ whenever $a \neq b$, $i_a = i_b$).

\medskip

\begin{example} 
\label{ex:finite type}

For $\fg$ of finite type, we have the following explicit description
\begin{equation}
	\label{eqn:e}
	\CS = \bigoplus_{\bn \in \nn} \CS_{\bn}, \text{ where } \CS_{(n_i \geq 0)_{i \in I}} = \Big\{ F(z_{i1},\dots,z_{in_i})_{i \in I} \text{ satisfying \eqref{eqn:wheel 1}} \Big\} 
\end{equation}
As $i \neq j$ range over $I$, the following are called Feigin-Odesskii wheel conditions:
\begin{equation}
	\label{eqn:wheel 1}
	F(\dots, z_{ia}, \dots,z_{jb},\dots)\Big|_{(z_{i1},z_{i2}, \dots, z_{i,1-c_{ij}}) \mapsto (z_{j1} q^{d_{ij}}, z_{j1} q^{d_{ij}+d_{ii}},  \dots, z_{j1} q^{-d_{ij}})} =  0
\end{equation}
The inclusion $\subseteq$ of \eqref{eqn:e} was established in \cite{E} based on the seminal work \cite{FO}, while the inclusion $\supseteq$ of \eqref{eqn:e} was proved in \cite{NT}. More generally, the formulas above hold as stated for strongly symmetrizable $\fg$, as shown in \cite[Section 3.25]{N new new}.

\end{example}

\medskip

\begin{example} 
\label{ex:simply laced}

For any simply-laced Kac-Moody $\fg$, it was shown in \cite{N Loop} that $\CS$ is the set of color-symmetric Laurent polynomials $F$ that satisfy the following condition: for any $i \neq j$ in $I$ and any $t-s, t'-s' \in 2 \BN$ such that $s-t' = s' -t \equiv d_{ij}$ mod 2, 
\begin{equation}
\label{eqn:wheel 2}
(x-y)^{\frac {t-s'+d_{ij}}2+1}
\end{equation}
divides the specialization of $F$ at 
\begin{align}
&z_{i1} = xq^s, \ \quad z_{i2} = xq^{s+2}, \quad \ \dots \ \quad z_{i,\frac {t-s}2} = xq^{t-2}, \ \ \quad z_{i, \frac {t-s}2+1} = xq^t \label{eqn:spec 1 intro} \\
&z_{j1} = yq^{s'}, \quad z_{j2} = yq^{s'+2}, \quad \dots \quad z_{j,\frac {t'-s'}2} = yq^{t'-2}, \quad z_{j, \frac {t'-s'}2+1} = yq^{t'} \label{eqn:spec 2 intro}
\end{align}

\end{example}

\medskip

\subsection{Simple modules and $q$-characters}
\label{sub:simple}

The isomorphism between $\UU$ and the double of the shuffle algebra $\CS$ allowed us to give the following explicit construction of an important class of simple modules (\cite{N Cat}). For any highest $\ell$-weight 
\begin{equation}
\label{eqn:bpsi}
\bpsi = \left( \psi_i(z) = \text{constant}_i \cdot \frac {\prod_{\gamma \in \BK^*} \left(1 - \frac {\gamma}z \right)^{l_{i,\gamma}}}{\prod_{\gamma \in \BK^*} \left(1 - \frac {\gamma}z \right)^{m_{i,\gamma}}} \right)_{i \in I}
\end{equation}
with various $l_{i,\gamma}, m_{i,\gamma} \geq 0$ (almost all of which are actually equal to 0) one considers the corresponding simple module $L(\bpsi)$ of the Borel half of $\UU$ (defined in \cite{HJ} in the finite type case and in \cite{N Cat} in the arbitrary symmetrizable Kac-Moody case). It was shown in \cite[Proposition 4.7]{N Char} that there is a isomorphism of vector spaces
\begin{equation}
\label{eqn:l}
L(\bpsi) \cong L^{\ord \bpsi} \otimes L^{\neq 0}(\bpsi)
\end{equation}
The first tensor factor in the RHS above was calculated in \cite{N Char} for finite type $\fg$ (proving a conjecture of \cite{MY, W}) and we will not be concerned with it. Instead, we will focus our attention on the second tensor factor in the RHS of \eqref{eqn:l}: this factor was realized in \cite[Definition 4.3, Remark 4.5]{N Char} as
\begin{equation}
\label{eqn:l neq 0}
L^{\neq 0}(\bpsi) = \bigoplus_{\bn \in \nn} \CS_{\bn} \Big/ \left(\CS_{\bn} \cap \Theta(\bpsi)_{\bn} \right)
\end{equation}
where $\Theta(\bpsi)_{\bn}$ denotes the set of $F(z_{i1}, \dots , z_{in_i})_{i \in I} \in \CV_{\bn}$ such that \footnote{Below, the symbol ``$\underset{z \in \BK^*}{\text{Res}}$" should be interpreted as ``the sum over all residues in $\BK^*$".}
\begin{equation}
\label{eqn:many residues}
\underset{z_n \in \BK^*}{\text{Res}} \dots \underset{z_1 \in \BK^*}{\text{Res}} \left(\frac {z_1^{d_1}\dots z_n^{d_n}F(z_1,\dots,z_n)}{\prod_{1\leq a < b \leq n} \zeta_{i_bi_a}\left(\frac {z_b}{z_a} \right)}  \prod_{a=1}^n \frac {\psi_{i_a}(z_a)}{z_a} \right) = 0
\end{equation}
for all orderings $i_1, \dots, i_n$ of $\bn$ and all $d_1,\dots,d_n \in \BZ$.  Given such an ordering, the notation $F(z_1,\dots,z_n)$ in the RHS of \eqref{eqn:many residues} refers to the fact that each symbol $z_a$ is matched with a variable of the form $z_{i_a \bullet_a}$ of $F$, where the choice of $\bullet_a$ does not matter due to the color-symmetry of $F$ (as long as $\bullet_a \neq \bullet_b$ whenever $a\neq b, i_a = i_b$).

\medskip

\noindent Since both $\CS_{\bn}$ and $\Theta(\bpsi)_{\bn}$ are preserved by multiplication with arbitrary elements of the ring $\CP_\bn^\pm$ of color-symmetric Laurent polynomials, we can decompose the right-hand side of \eqref{eqn:l neq 0} in terms of its fibers over the maximal spectrum
$$
\text{Spec}(\CP^\pm_{\bn}) \supseteq (\BK^*)^{\bn} := \prod_{i \in I} (\BK^*)^{n_i} / S_{n_i} 
$$
Explicitly, as $\bx = (x_{ia})_{i \in I, a \in \{1,\dots,n_i\}}$ runs over $(\BK^*)^{\bn}$, we have
\begin{equation}
\label{eqn:l neq 0 decompose}
L^{\neq 0}(\bpsi) = \mathop{\bigoplus_{\bn \in \nn}}_{\bx \in (\BK^*)^{\bn}} \left(\CS_{\bn} \Big/ \left(\CS_{\bn} \cap \Theta(\bpsi)_{\bn} \right) \right)_{\bx} = \mathop{\bigoplus_{\bn \in \nn}}_{\bx \in (\BK^*)^{\bn}} \CS_{\bn} \Big/ \left(\CS_{\bn} \cap \Theta(\bpsi)_{\bx} \right)
\end{equation}
where $\Theta(\bpsi)_{\bx}$ denotes the set of $F \in \CV_{\bn}$ such that 
\begin{equation}
\label{eqn:one residue}
\underset{z_n = x_n}{\text{Res}} \dots \underset{z_1 = x_1}{\text{Res}} \left(\frac {z_1^{d_1}\dots z_n^{d_n}F(z_1,\dots,z_n)}{\prod_{1\leq a < b \leq n} \zeta_{i_bi_a}\left(\frac {z_b}{z_a} \right)}  \prod_{a=1}^n \frac {\psi_{i_a}(z_a)}{z_a} \right) = 0
\end{equation}
for any ordering $x_1,\dots,x_n$ of (the entries of the tuple) $\bx$ underlain by any ordering $i_1,\dots,i_n$ of $\bn$ (in other words, each $x_a$ is a placeholder for one of the coordinates $x_{i_a \bullet_a}$ of $\bx$, where the numbers $\bullet_a$ are arbitrary as long as $\bullet_a \neq \bullet_b$ if $a \neq b$, $i_a = i_b$).

\medskip

\subsection{A model for shuffle algebras}
\label{sub:coefficients}

As explained in Subsection \ref{sub:shuffle intro}, we are interested in finding a geometric incarnation for the non-negative integers 
\begin{equation}
\label{eqn:the coefficient}
 \dim_{\BK} \left(\CS_{\bn} \Big/ \left(\CS_{\bn} \cap \Theta(\bpsi)_{\bx} \right) \right)
\end{equation}
for any (henceforth fixed) $\bn \in \nn$ and $\bx \in (\BK^*)^{\bn}$. We have the following inclusion
\begin{equation}
\label{eqn:inclusion}
\CS_{\bn} \Big/ \left(\CS_{\bn} \cap \Theta(\bpsi)_{\bx} \right) \hookrightarrow \CV_{\bn} \Big/  \Theta(\bpsi)_{\bx}
\end{equation}
and we will proceed to give a purely combinatorial description of the vector space in the right-hand side of \eqref{eqn:inclusion}. Define
\begin{equation}
\label{eqn:straight v}
V_{\bx} = \BQ [ z_{i1}^{\pm 1}, \dots, z_{in_i}^{\pm 1} ]_{i \in I}^{\text{sym}_{\bx}} 
\end{equation}
where $\sym_{\bx}$ refers to the partial symmetrization that only permutes those $z_{ia}$ and $z_{ib}$ for which $x_{ia} = x_{ib}$. For any tuple of integers $\bm = (m_{ia})_{i \in I, 1 \leq a \leq n_i}$, let
$$
\Theta_{\bx}^{\bm} \subseteq V_{\bx}
$$
be the intersection over all orderings $x_1,\dots,x_n$ of the entries of $\bx$ (with underlying ordering $i_1,\dots,i_n$ of $\bn$) of the sets of Laurent polynomials $F \in V_{\bx}$ such that 
\begin{multline}
\label{eqn:x residues}
\underset{z_n = 1}{\text{Res}} \dots \underset{z_1 = 1}{\text{Res}} \left(\frac {z_1^{d_1}\dots z_n^{d_n}F(z_1,\dots,z_n)}{\prod_{a=1}^n (z_a - 1)^{m_{a}}} \right. \\ \left. \prod_{1\leq a < b \leq n} (z_b - z_a)^{\delta_{(i_a,x_a) = (i_b,x_b)} - \delta_{x_a = x_b q^{d_{i_ai_b}}}} \right) = 0
\end{multline}
for all $d_1,\dots,d_n \in \BZ$, where the numbers $m_1,\dots,m_n$ represent an ordering of the entries of $\bm$ in the same way as $x_1,\dots,x_n$ represent an ordering of the entries of $\bx$.

\medskip

\begin{remark}
\label{rem:the quiver}

If $m_{ia} \geq 0, \forall i,a$, then $\Theta_{\bx}^{\bm}$ coincides with \eqref{eqn:kernel} for the quiver with 

\medskip

\begin{itemize}[leftmargin=*]

\item vertex set $I_{\bx} = \left\{ (i,a) \Big| i \in I, 1 \leq a \leq n_i \right\} \Big / \Big((i,a) \sim (i,b) \text{ if } x_{ia} = x_{ib} \Big)$;

\medskip

\item arrow set $(i,a) \rightarrow (j,b)$ whenever $x_{ia} = x_{jb} q^{d_{ij}}$;

\medskip

\item dimension vector $(i,a) \leadsto$ number of times $x_{ia}$ appears in the multiset $\{x_{i1},\dots,x_{in_i}\}$ 

\medskip

\item framing dimension $m_{ia}$ at any vertex $(i,a)$.

\end{itemize}

\end{remark}

\medskip

\begin{proposition}
\label{prop:rewrite}

We have an isomorphism of finite-dimensional vector spaces
\begin{equation}
\label{eqn:rewrite}
\CV_{\bn} \Big/ \Theta(\bpsi)_{\bx} \cong \left( V_{\bx} \Big/ \Theta^{\bm}_{\bx} \right) \otimes_{\BQ} \BK
\end{equation}
where the tuple $\bm$ is related to $\bpsi$ by the formula $\psi_i(z) = (z-x_{ia})^{-m_{ia}}$ times a rational function which is regular and non-zero at $z = x_{ia}$.

\end{proposition}

\medskip

\begin{proof} There is a natural map
\begin{equation}
\label{eqn:injection}
\CV_{\bn} \hookrightarrow V_{\bx}  \otimes_{\BQ} \BK
\end{equation}
defined by rescaling the variables of Laurent polynomials according to
$$
z_{ia} \leadsto z_{ia}x_{ia}
$$
and then regarding color-symmetric Laurent polynomials as being simply partially symmetric with respect to those $z_{ia}$ and $z_{ib}$ such that $x_{ia} = x_{ib}$. To show that the assignment \eqref{eqn:injection} sends $\Theta(\bpsi)_{\bx}$ to $\Theta_{\bx}^{\bm}$, we will argue as follows: a Laurent polynomial $F$ lies inside $\Theta(\bpsi)_{\bx}$ if and only if finitely many linear combinations of its derivatives vanish at $\bx$. This vanishing property is preserved upon multiplying $F$ by any rational function which is regular and non-zero at $\bx$. However, this is precisely the relation between the integrands of \eqref{eqn:one residue} and \eqref{eqn:x residues}: the latter expression involves the linear factors of the former expression which vanish at the point $\bx$ and nothing more. Thus, $F$ lies in $\Theta(\bpsi)_{\bx}$ if and only if its image under \eqref{eqn:injection} lies in $\Theta_{\bx}^{\bm}$.

\medskip 

\noindent The previous argument precisely shows that \eqref{eqn:injection} gives rise to an injective map
\begin{equation}
\label{eqn:injection quotients}
\CV_{\bn} \Big/ \Theta(\bpsi)_{\bx} \hookrightarrow \left( V_{\bx} \Big/ \Theta^{\bm}_{\bx} \right) \otimes_{\BQ} \BK
\end{equation}
We will show that \eqref{eqn:injection quotients} is also surjective by the following general argument. Consider any parabolic subgroup (in the case at hand, $P$ will be the stabilizer of $\bx$)
$$
P = S_{n_1} \times \dots \times S_{n_k} \hookrightarrow S_n
$$
where $n = n_1+\dots+n_k$, and any point
$$
x = (x_1,\dots,x_n) \in (\BC^*)^n
$$
such that $x_i = x_j$ if and only if $(ij) \in P$. Then for any $P$-invariant Laurent polynomial $f(z_1,\dots,z_n)$, we \underline{claim} that there exists a $S_n$-invariant Laurent polynomial $g$ which has the same $N$-th order (for arbitrarily large $N$) derivatives at $x$ as $f$; it is clear that the underlined claim precisely establishes the surjectivity of \eqref{eqn:injection quotients}. To prove the underlined claim, it suffices to take
$$
g = \sum_{\sigma \in S_n / S_{n_1} \times \dots \times S_{n_k}} \sigma(f\cdot r)
$$
where $r(z_1,\dots,z_n)$ is a $P$-invariant Laurent polynomial such that
$$
r \in \Big(1 + \left( (z_1 - x_1)^N, \dots, (z_n - x_n)^N \right) \Big) \bigcap_{(ij) \notin P} \Big( (z_i - x_j)^N \Big)
$$
The existence of such a Laurent polynomial $r$ is due to the fact that the polynomial
$$
p(z_1,\dots,z_n) = \prod_{(ij) \notin P} (z_i - x_j)^N 
$$
does not vanish at $(z_1,\dots,z_n) = (x_1,\dots,x_n)$, and thus the principal ideal $(p)$ together with $( (z_1 - x_1)^N, \dots, (z_n - x_n)^N)$ generate the unit ideal.  \end{proof}

\bigskip

\section{Relation to critical $K$-theory}
\label{sec:k-ha}

\medskip

\noindent In Subsections \ref{sub:shuffle intro} and \ref{sub:simple}, we saw that the coefficients of $q$-characters of simple modules $L(\bpsi)$ can be obtained from the dimensions of the vector spaces
\begin{equation}
\label{eqn:vector k-ha}
\CS_{\bn} \Big/ \left(\CS_{\bn} \cap \Theta(\bpsi)_{\bx} \right)
\end{equation}
as $\bn$ ranges over $\nn$ and $\bx$ ranges over $(\BK^*)^{\bn}$. We will now relate the vector spaces \eqref{eqn:vector k-ha} to critical $K$-theory of quiver varieties, which we will argue provides a suitable generalization of Theorem \ref{thm:hl} beyond reachable simple modules. We assume that $\fg$ is strongly symmetrizable throughout the present Section; recall that this includes all finite and affine type symmetrizable Kac-Moody Lie algebras except for $\widehat{\fsl}_2$.

\medskip

\subsection{Quiver with potential} 
\label{sub:quiver potential}

Let us recall from Subsection \ref{sub:k-ha intro} the quiver $Q$ with

\medskip

\begin{itemize}[leftmargin=*]

\item vertex set $I$,

\medskip 

\item arrows $j \xleftarrow{_j\square_{i}} i$ for all $i,j \in I$ (including $i=j$),

\medskip

\item potential (i.e. a linear combination of cycles in the path algebra of $Q$)
\begin{equation}
\label{eqn:potential}
W = \sum_{i \neq j} {_j\square_{i}}  \underbrace{ {_i\square_{i}} \dots {_i\square_{i}}}_{-c_{ij} \text{ factors}} {_i\square_{j}}
\end{equation}
where $c_{ij} = \frac {2d_{ij}}{d_{ii}} \in \BZ_{\leq 0}$ are the off-diagonal coefficients of the Cartan matrix.

\end{itemize}

\medskip

\noindent In \cite{VV1,VV2}, the authors considered the critical $K$-theory of the above quiver with potential (following \cite{BFK, EP, Hi}, see also \cite{Pad}). We refer to the aforementioned works for an introduction to critical $K$-theory, and restrict ourselves to simply listing the following key facts that will be used in the present paper.

\medskip

\begin{itemize}[leftmargin=*]

\item For any $\bn \in \nn$, recall the stack $M_{\bn}$ of $\bn$-dimensional representations of $Q$. There is an action $\BC^* \curvearrowright M_{\bn}$ which assigns weight $q^{d_{ij}}$ to any arrow $_j\square_i$.

\medskip

\item Evaluating the trace of the potential $W$ in quiver representations gives a function
\begin{equation}
\label{eqn:function}
\text{Tr}(W) : M_{\bn} \rightarrow \BC
\end{equation}
which is $\BC^*$-invariant. As such, there is  an equivariant critical $K$-theory group
\begin{equation}
\label{eqn:critical k-theory}
K^{\BC^*}(Q,W)
\end{equation}
All $K$-theory groups (ordinary or critical) in the present paper will be considered to be localized, i.e. $\otimes_{\BZ} \BQ$ in the non-equivariant case and $\otimes_{\BZ[q^{\pm 1}]} \BQ(q)$ in the $\BC^*$-equivariant case. Thus, the abelian group \eqref{eqn:critical k-theory} is a $\BK$-vector space.

\medskip

\item The critical $K$-theory \eqref{eqn:critical k-theory} is a module over the equivariant $K$-theory of the stack
$$
K^{\BC^*}(M_{\bn}) = \BK[z_{i1}^{\pm 1}, \dots, z_{in_i}^{\pm 1}]_{i \in I}^{\text{sym}} = \CV_{\bn}
$$

\medskip

\item There is a $\CV_{\bn}$-module homomorphism
\begin{equation}
\label{eqn:map to shuffle}
K^{\BC^*}(Q,W)_{\bn} \xrightarrow{\iota_{\bn}} \CV_{\bn}
\end{equation}
(\cite[Proposition 3.6]{Pad Cat}) whose image is contained in the image of the push-forward
\begin{equation}
\label{eqn:push forward}
K^{\BC^*}(Z_{\bn}) \xrightarrow{\alpha} K^{\BC^*}(M_{\bn}) 
\end{equation}
with $Z_{\bn}$ denoting the zero locus of the function $\text{Tr}(W) : M_{\bn} \rightarrow \BC$.

\medskip

\item There is a convolution algebra structure (\cite{VV1, VV2}, see also \cite{Pad}) on
\begin{equation}
\label{eqn:item 3}
K^{\BC^*}(Q,W) = \bigoplus_{\bn \in \nn} K^{\BC^*}(Q,W)_{\bn}
\end{equation}
that is referred to as the critical convolution algebra. The map
\begin{equation}
\label{eqn:item 4}
K^{\BC^*}(Q,W)  \xrightarrow{\iota} \CV
\end{equation}
obtained by putting together all the maps \eqref{eqn:map to shuffle} is an algebra homomorphism, where $\CV$ is made into an algebra with respect to the shuffle product \eqref{eqn:mult}. 

\end{itemize}

\medskip

\noindent The above notions are particularly simple when $\bn = \bs^i$ for some $i \in I$. In this case, the moduli stack of quiver representations is $\BC \times \cdot /\BC^*$, and so its $K$-theory is 
$$
K^{\BC^*}(Q,W)_{\bs^i} = K^{\BC^*}(M_{\bs^i}) \cong \BK[z_{i1}^{\pm 1}] = \CV_{\bs^i}
$$
In particular, since the $\BK$-algebra $\CS$ is generated by $\{z_{i1}^d\}_{i \in I, d \in \BZ}$, we have
\begin{equation}
\label{eqn:one way}
\text{Im } \iota \supseteq \CS
\end{equation}

\medskip

\begin{proposition} \label{prop:yu} If $\fg$ is strongly symmetrizable, then $\emph{Im }\iota = \CS$. \end{proposition}

\medskip

\begin{proof} We follow the idea of Yu Zhao from \cite{Zhao}. Recall the explicit description of $\CS$ from Example \ref{ex:finite type}. It suffices to prove that $\text{Im }\alpha \subseteq \CS_{\bn}$ with $\alpha$ being the map \eqref{eqn:push forward}. Because of the localization exact sequence
$$
K^{\BC^*}(Z_{\bn}) \xrightarrow{\alpha} K^{\BC^*}(M_{\bn}) \cong \CV_{\bn} \xrightarrow{\beta} K^{\BC^*}(M_{\bn} \backslash Z_{\bn})
$$
we must prove that $\text{Ker }\beta \subseteq \CS_{\bn}$. To this end, fix $i \neq j$ in $I$ and consider the locally closed substack $U \subset M_{\bn}$ that parameterizes quiver representations of the form
\begin{equation}
\label{eqn:triangle}
\xymatrix{& & \BC \ar[rrd]^{_i\square_j} & & \\
\BC \ar[rru] ^{_j\square_i} & \BC \ar[l]^{_i\square_i} & \dots \ar[l]^{_i\square_i} &  \BC \ar[l]^{_i\square_i} & \BC \ar[l]^{_i\square_i}}
\end{equation}
(the top $\BC$ lies above the vertex $j$ and the direct sum of the bottom $\BC$'s lies above the vertex $i$) where the arrows displayed are all $\neq 0$, direct sum any collection of vector spaces such that all arrows to/from them are 0. Our choice of arrows in \eqref{eqn:triangle} implies that $U \subset M_{\bn} \backslash Z_{\bn}$, and thus it suffices to show that any Laurent polynomial whose restriction to $U$ is 0 actually lies in $\CS_{\bn}$. However, 
$$
U \cong (\BC^*)^{2-c_{ij}} \Big / (\BC^*)^{2-c_{ij}} \cong \BC^* \times \left( \cdot / \BC^* \right)
$$ 
where the left-most $(\BC^*)^{2-c_{ij}}$ parameterizes the arrow maps in \eqref{eqn:triangle} and the right-most $(\BC^*)^{2-c_{ij}}$ acts on it by conjugation. Restricting a Laurent polynomial in $\CV_{\bn} \cong K^{\BC^*}(M_{\bn})$ to $U$ has the effect of performing the specialization \eqref{eqn:wheel 1}.

\end{proof}

\medskip

\subsection{Graded quiver with potential}
\label{sub:graded quiver potential}

Let us now consider the graded quiver $Q^{\gr}$ of Subsection \ref{sub:euler intro}, which has

\medskip

\begin{itemize}[leftmargin=*]

\item vertex set $I \times \BZ$,

\medskip 

\item arrows $(j,t) \xleftarrow{^{t}_j\square^{s}_{i}} (i,s)$ for all $i,j \in I$ and $s,t \in \BZ$ with $t = s-d_{ij}$;

\medskip

\item potential
$$
W^{\gr} = \sum_{i \neq j} \sum_{s \in \BZ} \left( ^{s}_j\square_{i}^{s+d_{ij}} \right) \left(^{s+d_{ij}}_i\square_{i}^{s+d_{ij}+d_{ii}} \right) \dots \left( ^{s-d_{ij}-d_{ii}}_i\square_{i}^{s-d_{ij}} \right) \left(^{s-d_{ij}}_i\square_{j}^{s} \right)
$$

\end{itemize}

\medskip

\noindent The connection between the graded quiver above and the construction of Subsection \ref{sub:quiver potential} is given via the fixed points of the action $\BC^* \curvearrowright M_{\bn}$. More specifically, we have
\begin{equation}
\label{eqn:fixed points stack}
M_{\bn}^{\BC^*} \supset \bigsqcup_{\bx \in q^{\BZ^{\bn}}} M_{\bx} 
\end{equation}
where in the RHS, $M_{\bx}$ denotes the stack of $Q^{\text{gr}}$-representations which have dimension at $(i,s)$ equal to the number of times $q^s$ appears among the coordinates of $\bx = (x_{ia})_{i \in I, 1 \leq a \leq n_i}$. The map \eqref{eqn:injection} can be interpreted as the restriction from the $K$-theory of $M_{\bn}$ to that of $M_{\bx}$ (note that $K(M_{\bx}) = V_{\bx}$). We have the following non-equivariant analogue of the objects in \eqref{eqn:map to shuffle}
\begin{equation}
\label{eqn:iota gr}
 K(Q^{\gr}, W^{\gr})_{\bx} \xrightarrow{\iota^{\gr}_{\bx}} V_{\bx}
\end{equation}

\medskip

\subsection{Framed critical $K$-theory}
\label{sub:framed k-ha}

Let us henceforth consider a highest $\ell$-weight $\bpsi$ of the form
$$
\bpsi = \left( \psi_i(z) = \text{constant}_i \cdot \frac {\psi_i^{\text{num}}(z)}{\psi_i^{\text{den}}(z)} \right)_{i \in I}
$$
as in \eqref{eqn:bpsi intro}. In other words, we have that $\bpsi^{\text{num}} = (\psi_i^{\text{num}}(z))_{i \in I}$, $\bpsi^{\text{den}} = (\psi_i^{\text{den}}(z))_{i \in I}$ are Laurent polynomial $\ell$-weights of the form
\begin{align*}
&\psi_i^{\text{num}}(z) = \prod_{s \in \BZ} \left(1- \frac {q^s}{z} \right)^{l_{i,s}} \\
&\psi_i^{\text{den}}(z) = \prod_{s \in \BZ}  \left(1- \frac {q^s}{z} \right)^{m_{i,s}}
\end{align*}
for various $l_{i,s}, m_{i,s} \in \BN$, almost all of which are 0. We will write $\bl = (l_{i,s})_{i \in I, s \in \BZ}$ and $\bm = (m_{i,s})_{i \in I, s \in \BZ}$. Consider the quivers from the previous Subsections, but framed as follows:

\medskip

\begin{itemize}[leftmargin=*]
	
\item we endow $Q$ with framing $\oplus_{s \in \BZ} \BC^{m_{i,s}}$ at every vertex $i \in I$, and declare the $s$-th direct summand to be equivariant with respect to the character $q^s$ of $\BC^*$.

\medskip

\item we endow $Q^{\gr}$ with framing $\BC^{m_{i,s}}$ at every vertex $(i,s) \in I \times \BZ$.

\end{itemize}

\medskip 

\noindent The choices above allow to define the stacks of framed quiver representations $M_{\bn,\bm}$ and $M_{\bx,\bm}$, which are total spaces of affine bundles over the stacks $M_{\bn}$ and $M_{\bx}$, respectively. As such, the corresponding framed critical $K$-theory groups are identified with the unframed ones, i.e. in the following equations 
\begin{equation}
\label{eqn:map to shuffle framed}
K^{\BC^*}(Q,W)_{\bn,\bm} \xrightarrow{\iota_{\bn,\bm}} \CV_{\bn}
\end{equation}
\begin{equation}
\label{eqn:iota gr framed}
K(Q^{\gr}, W^{\gr})_{\bx,\bm} \xrightarrow{\iota^{\gr}_{\bx,\bm}} V_{\bx}
\end{equation}
both the LHS and the map are identified with those of \eqref{eqn:map to shuffle} and \eqref{eqn:iota gr}, respectively.

\medskip

\noindent The Laurent polynomial $\ell$-weight $\bpsi^{\text{num}}$ gives rise to the following element 
\begin{equation}
\label{eqn:delta}
\Delta_{\bl} = \prod_{i \in I} \prod_{s \in \BZ} \wedge^\bullet \left(\frac {q^s}{\CT_{i}} \right)^{l_{i,s}} 
\end{equation}
in the $K$-theory of $M_{\bn,\bm}$ (a.k.a. $\CV_{\bn}$), where $\CT_{i}$ denotes the tautological vector bundle corresponding to the vertex $i \in I$. Similarly, we have the element
\begin{equation}
\label{eqn:delta graded}
\Delta^{\gr}_{\bl} = \prod_{i \in I} \prod_{s \in \BZ} \wedge^\bullet \left(\frac 1{T_{i,s}} \right)^{l_{i,s}}
\end{equation}
in the $K$-theory of $M_{\bx,\bm}$ (a.k.a. $V_{\bx}$), where $T_{i,s}$ denotes the tautological vector bundle corresponding to the vertex $(i,s) \in I \times \BZ$. The relation between $\Delta_{\bl}$ and $\Delta_{\bl}^{\gr}$ is that the latter is the product of those linear factors of the former which do not involve any $q$'s, upon invoking the restriction map
$$
[\CT_i] \Big|_{\BC^* \text{ fixed locus}} = \sum_{s \in \BZ} q^{s} [T_{i,s}]
$$

\medskip

\begin{lemma}
\label{lem:multiplication}

Under the maps \eqref{eqn:map to shuffle framed} and \eqref{eqn:iota gr framed}, the operations
$$
K^{\BC^*}(Q,W)_{\bn,\bm} \xrightarrow{\cdot \Delta_{\bl}} K^{\BC^*}(Q,W)_{\bn,\bm} 
$$
$$
K(Q^{\egr}, W^{\egr})_{\bx,\bm} \xrightarrow{\cdot \Delta_{\bl}^{\egr}} K(Q^{\egr}, W^{\egr})_{\bx,\bm}
$$
correspond to multiplication by elements which we abusively denote by
\begin{equation}
\label{eqn:delta 1}
\Delta_{\bl} = \prod_{i \in I} \prod_{a=1}^{n_i} \prod_{s \in \BZ} \left(1 - \frac {q^s}{z_{ia}} \right)^{l_{i,s}} \in \CV_{\bn} 
\end{equation} 
\begin{equation}
\label{eqn:delta 2}
\Delta^{\egr}_{\bl} = \prod_{i \in I} \prod_{a=1}^{n_i} \left(1 - \frac {1}{z_{ia}} \right)^{l_{i,\log_q(x_{ia})}} \in V_{\bx} 
\end{equation}
respectively. The relation between $\Delta_{\bl}$ and $\Delta_{\bl}^{\egr}$ is that the latter consists of all the linear factors of the former which vanish at $\{z_{ia} = x_{ia}\}_{i \in I, 1 \leq a \leq n_i}$ \footnote{The apparent mismatch in notation between $\Delta_{\bl}$ and $\Delta_{\bl}^{\gr}$ is due to the map $\CV_{\bn} \rightarrow V_{\bx} \otimes_{\BQ} \BK$ of \eqref{eqn:injection} shifting the variables according to $z_{ia} \leadsto z_{ia} x_{ia}$.}.

\end{lemma}

\medskip

\subsection{Framed stable critical $K$-theory} 
\label{sub:framed stable k-ha}

Following a suggestion of Shivang Jindal, we restrict the objects \eqref{eqn:map to shuffle framed} and \eqref{eqn:iota gr framed} to the stable open subsets $M_{\bn,\bm}^{\stab} \subset M_{\bn,\bm}$ and $M_{\bx,\bm}^{\stab} \subset M_{\bx,\bm}$, thus obtaining the framed stable critical $K$-theory groups \footnote{Stable framed critical $K$-theory was studied in \cite{VV1,VV2}, see also \cite{DM, S} in cohomology.}
\begin{equation}
\label{eqn:open locus}
K^{\BC^*}(Q, W)^{\stab}_{\bn,\bm} \xrightarrow{\iota^{\stab}_{\bn,\bm}} K^{\BC^*}(M_{\bn,\bm}^{\stab}) \cong \CV_{\bn} \Big / \Theta \left(\frac 1{\bpsi^{\text{den}}} \right)_{\bn} 
\end{equation}
\begin{equation}
\label{eqn:open locus gr}
K(Q^{\gr}, W^{\gr})^{\stab}_{\bx,\bm} \xrightarrow{\iota^{\gr,\stab}_{\bx,\bm}} K(M_{\bx,\bm}^{\stab}) \cong V_{\bx} \Big/ \Theta_{\bx}^{\bm}
\end{equation}
The right-most isomorphism in \eqref{eqn:open locus gr} is simply Theorem \ref{thm:nilp}. Meanwhile, the right-most isomorphism in \eqref{eqn:open locus} can either be obtained by developing an equivariant analogue of Theorem \ref{thm:nilp}, or by deducing it from the right-most isomorphism in \eqref{eqn:open locus gr} via equivariant localization and Proposition \ref{prop:rewrite}. We thus obtain the following diagram, where the top vertical arrows denote restriction to the open locus, and the bottom vertical arrows denote restrictions to the $\BC^*$-fixed locus:
\begin{equation}
\label{eqn:big diagram}
\xymatrix{
K^{\BC^*}(Q,W)_{\bn} \ar@/_3pc/[dd]_{\pi} \ar[d] \ar[r]^-{\cdot \Delta_{\bl}} & K^{\BC^*}(Q,W)_{\bn} \ar[d] \ar[r]^-{\iota_{\bn}} & \CV_{\bn} \ar@{->>}[d]  \\
K^{\BC^*}(Q,W)_{\bn}^{\stab} \ar[d] \ar[r]^-{\cdot \Delta_{\bl}} & K^{\BC^*}(Q,W)^{\stab}_{\bn} \ar[d] \ar[r]^-{\iota_{\bn}^{\stab}} & \CV_{\bn} /  \Theta \left(\frac 1{\bpsi^{\text{den}}} \right)_{\bn}  \ar[d]_{\text{Prop. \ref{prop:rewrite}}} \\
K(Q^{\gr}, W^{\gr})^{\stab}_{\bx,\bm} \otimes_{\BQ} \BK \ar[r]^-{\cdot \Delta_{\bl}^{\gr}} & K(Q^{\gr}, W^{\gr})^{\stab}_{\bx,\bm} \otimes_{\BQ} \BK \ar[r]^-{\iota_{\bx,\bm}^{\gr,\stab}} & V_{\bx} / \Theta_{\bx}^{\bm} \otimes_{\BQ}  \BK}
\end{equation}
The top squares commute, while the bottom ones do not.

\medskip

\begin{proof} \emph{of Theorem \ref{thm:main intro}:} We showed in Proposition \ref{prop:yu} that $\text{Im }\iota_{\bn} = \CS_{\bn}$. Thus, the image of the composition of the top two horizontal maps is the vector space $\Delta_{\bl} \cdot \CS_{\bn}$. Mapping further down the two vertical maps, we obtain a vector space isomorphic to
$$
\Delta_{\bl} \cdot \CS_{\bn} \Big / \left( \Delta_{\bl} \cdot \CS_{\bn} \cap \Theta \left(\frac 1{\bpsi^{\text{den}}} \right)_{\bx} \right) \ \ \subset \ \ \CV_{\bn} \Big/ \Theta \left(\frac 1{\bpsi^{\text{den}}} \right)_{\bx}  \stackrel{\text{Prop. \ref{prop:rewrite}}}{\cong} V_{\bx} / \Theta_{\bx}^{\bm} \otimes_{\BQ}  \BK
$$
As in \eqref{eqn:2}, multiplication by $\Delta_{\bl}$ gives an isomorphism between the LHS above and
$$
\CS_{\bn} \Big / \left(\CS_{\bn} \cap \Theta \left(\frac {\bpsi^{\text{num}}}{\bpsi^{\text{den}}} \right)_{\bx} \right)
$$
precisely as we needed to show. \end{proof}

\medskip

\subsection{Coda}
\label{sub:conj}

We have just shown that we can obtain a vector space naturally isomorphic to \eqref{eqn:main space} as the image of the composition $\rightarrow \rightarrow \downarrow \downarrow$ in diagram \eqref{eqn:big diagram}. 

\medskip

\begin{lemma}
\label{lem:have the same image}

The compositions $\rightarrow \rightarrow \downarrow\downarrow$ and $\downarrow \downarrow\rightarrow\rightarrow$ in \eqref{eqn:big diagram} have the same image.

\end{lemma}

\medskip

\begin{proof} Since the top squares in \eqref{eqn:big diagram} are commutative, it suffices to focus our attention on the bottom squares, in which the vertical arrows are given by restriction to the $\BC^*$-fixed locus and the horizontal arrows are given by push-forward and multiplication by $\Delta_{\bl}, \Delta_{\bl}^{\gr}$. The equivariant localization theorem  (see \cite[Proposition 2.6.(d)]{VV2} for the critical $K$-theory version) implies that the two compositions
$$
\Big[ \text{middle } K^{\BC^*}(Q,W)^{\stab}_{\bn} \Big]\rightarrow V_{\bx} / \Theta_{\bx}^{\bm} \otimes_{\BQ}  \BK
$$
in diagram \eqref{eqn:big diagram} have the same image, because restriction to the $\BC^*$-fixed locus is the inverse of the push-forward of the $\BC^*$-fixed locus once we invert non-trivial tangent $\BC^*$-weights (non-trivial tangent weights are invertible in the ring $V_{\bx} / \Theta_{\bx}^{\bm} \otimes_{\BQ}  \BK$, which is finite-dimensional and supported at $\{z_{ia} = 1\}_{i \in I, 1 \leq a \leq n_i}$). For a similar reason, because
$$
\Delta_{\bl} \Big|_{\BC^* \text{ fixed locus}} = \Delta_{\bl}^{\gr} \cdot u
$$
where $u$ is invertible in the ring $V_{\bx} / \Theta_{\bx}^{\bm} \otimes_{\BQ}  \BK$, all compositions
$$
\Big[ \text{leftmost } K^{\BC^*}(Q,W)^{\stab}_{\bn} \Big] \rightarrow V_{\bx} / \Theta_{\bx}^{\bm} \otimes_{\BQ}  \BK
$$
in diagram \eqref{eqn:big diagram} have the same image. 

\end{proof}

\medskip

\noindent Lemma \ref{lem:have the same image} implies that we may realize the vector space \eqref{eqn:main space} as the image of the composition $\downarrow \downarrow \rightarrow \rightarrow$ in \eqref{eqn:big diagram}, and we will argue that this composition is more natural than $\rightarrow \rightarrow \downarrow\downarrow$ from a geometric point of view. Indeed, it was proved in \cite{COZZ} that 
\begin{equation}
\label{eqn:surj inj}
\pi \text{ is surjective and } \iota_{\bx,\bm}^{\gr,\stab} \text{ is injective}
\end{equation}
using their theory of critical stable envelopes. Therefore, the image of the composition $\downarrow \downarrow \rightarrow \rightarrow$ in \eqref{eqn:big diagram} is isomorphic to the image of the map $\cdot \Delta_{\bl}^{\gr}$ therein. Thus
\begin{equation}
\label{eqn:complete resolution}
\CS_{\bn} \Big/ \left(\CS_{\bn} \cap \Theta(\bpsi)_{\bx}\right) \cong \left[ \Delta^{\gr}_{\bl} \cdot K(Q^{\gr}, W^{\gr})^{\stab}_{\bm,\bx} \right] \otimes_{\BQ} \BK
\end{equation}
which is the content of Corollary \ref{cor:main}. Since the left-hand side of the above equation categorifies the coefficients of the $q$-character, \eqref{eqn:complete resolution} is precisely the categorification and generalization to all $\bpsi$ of Theorem \ref{thm:hl}: since critical $K$-theory of quivers with potential is a natural replacement of the $K$-theory of the critical locus, then 
$$
K(Q^{\gr}, W^{\gr})^{\stab}_{\bx,\bm} \quad \text{is a replacement of} \quad K(M_{\bx,\bm}^{\stab} \cap \text{Crit}(\text{Tr}(W^{\gr})))
$$
Meanwhile, multiplication by $\Delta_{\bl}^{\gr}$ is a derived version of imposing the vanishing of $l_{i,s}$ generic linear combinations of paths ending at any vertex $(i,s)$. Therefore, the RHS of \eqref{eqn:complete resolution} is a natural replacement of $K(N_{\bx,\bm,\bl}^{\stab}) \otimes_{\BQ} \BK$.

\end{document}